\documentclass[a4paper,12pt]{article}
\usepackage[cp1251]{inputenc}
\usepackage{amssymb,amsmath,amsthm}
\usepackage[russian]{babel}
\usepackage{indentfirst}
\usepackage{amssymb,amsmath,amsthm}
\usepackage{hyperref,hypbmsec}
\usepackage[pdftex]{graphicx}

\numberwithin{equation}{section}

\newtheorem{Th}{\hskip\parindent Theorem}[section]
\newtheorem{Le}{\hskip\parindent Lemma}[section]

\newtheorem{Zam}{\hskip\parindent Remark}[section]
\newtheorem{Hyp}{\hskip\parindent Conjecture}[section]

\newcommand{\A}{\mathcal{A}}
\newcommand{\E}{\mathfrak{C}}
\newcommand{\R}{\mathfrak{R}}

\newcommand{\D}{\mathfrak{D}}
\newcommand{\N}{\mathbb{N}}
\newcommand{\M}{\mathfrak{M}}
\newcommand{\sys}{\mathfrak{S}}

\newcommand{\1}{\mathbf{1}}

\newcommand{\Om}{\widetilde{\Omega}}
\newcommand{\q}{\mathbf{q}}
\newcommand{\KK}{\overline{K}}

\newcommand{\rr}{\mathbb{R}}

\newcounter{propet}

\renewcommand{\le}{\leqslant}\renewcommand{\ge}{\geqslant}
\renewcommand{\proofname}{Док}

\begin{document}

\author{ D.\,A.\,Frolenkov\footnote{Research is supported by RFFI (grant № 11-01-00759-а)} \quad
I.\,D.\,Kan,\footnote{Research is supported by RFFI  РФФИ (grant № 12-01-00681-а)} }
\title{
\begin{flushright}
\small{Dedicated to the memory of\\professor N.M. Korobov.}
\end{flushright}
A reinforcement of the Bourgain-Kontorovich's theorem by elementary methods.
}
\date{}
\maketitle
УДК
\begin{abstract}
Recently (in 2011) several new theorems concerning this conjecture were proved by Bourgain and Kontorovich. The easiest of them states that the set of numbers satisfying Zaremba's conjecture with $A=50$ has positive proportion in $\N.$ The proof of this theorem is rather complicated and refers to the spectral theory. In this paper,using only elementary methods, the same theorem is proved with $A=13$ .

Bibliography: 17 titles.

 \par
\end{abstract}

\renewcommand{\proofname}{{\bf Доказательство}}

\setcounter{Zam}0
\section{Introduction}
Let $\R_{A}$ be the set of rational numbers whose continued fraction expansion has all partial quotients being bounded by  $A.$
Let $\D_{A}$ be the set of denominators of numbers in $\R_{A}$:
$\D_{A}=\left\{d\Bigl| \exists b:\, (b,d)=1, \frac{b}{d}\in\R_{A}\right\}.$
And set
$$
\D_{A}(N)=\left\{d\in\D_{A}\Bigl| d\le N \right\}.
$$
\begin{Hyp}(\textbf{Zaremba's conjecture} ~\cite[p. 76]{Zaremba},\,1971 ).
For sufficiently large $A$  one has
$$\D_{A}=\N.$$
\end{Hyp}
Let $\A \in \N$ be any finite alphabet ($|\A|\ge2$) and let $\R_{\A}$ and $\E_{\A}$  be the set of finite and infinite continued fractions whose partial quotients belong to $\A.$ And let
$$\D_{\A}(N)=\left\{d\Bigl| d\le N,\, \exists b: (b,d)=1,\, \frac{b}{d}\in\R_{\A}\right\}$$
be the set of denominators bounded by $N.$ Let $\delta_{\A}$ be the Hausdorff dimension of the set $\E_{\A}$ (another definition of $\delta_{\A}$ will be given below in §\ref{good function and kusic th}). In the article ~\cite{FK}, using the method, developed by Bourgain-Kontorovich in~\cite{BK}, as the base the following theorem was proved
\begin{Th}\label{uslov}
For any alphabet $\A$ with
\begin{equation}\label{KFcondition1}
\delta_{\A}>1-\frac{27-\sqrt{633}}{16}=0,8849\ldots,
\end{equation}
the following inequality (positive proportion) holds
\begin{equation}\label{KFresult}
\#\D_{\A}(N)\gg N.
\end{equation}
\end{Th}
\begin{Zam}
It is proved ~\cite{Jenkinson} that $\delta_{7}=0,8889\ldots.$ From this follows that the alphabet $\left\{1,2,\ldots,7\right\}$
satisfies the condition of Theorem \ref{uslov}. It is also proved in~\cite{Jenkinson} that for the alphabet $\A=\left\{1,2,\ldots,6,8\right\}$ one has $\delta_{\A}=0,8851\ldots.$ Consequently, the alphabet $\A=\left\{1,2,\ldots,6,8\right\}$ also satisfies the condition of Theorem \ref{uslov}.
\end{Zam}
Note that the Theorem \ref{uslov} is based on the Lemma 7.1 in the paper ~\cite{BK}, which we are going to state next. To begin with we describe all necessary objects.\par
Let $K, X, Y \ge1$ be real numbers and $q$ be a positive integer.
Moreover, let $\eta=(x,y)^t, \eta'=(u,v)^t \in \mathbb{Z}^2$ be vectors such as
$$|\eta|\asymp\frac{X}{Y}, |\eta'|\asymp X, (x,y)=1,(u,v)=1.$$
\begin{Le}\label{BK lemma 7.1}(~\cite[p. 46, Lemma 7.1.]{BK})
If the following inequality
\begin{equation}\label{BK lemma 7.1 condition}
(qK)^{\frac{13}{5}}<Y<X,
\end{equation}
holds, then
\begin{equation}\label{BK lemma 7.1 statement}
\#\left\{\gamma\in SL_2(\mathbb{Z})\Bigl| \|\gamma\|\asymp Y, |\gamma\eta-\eta'|<\frac{X}{K}, \gamma\eta\equiv\eta'\pmod{q}\right\}\ll\frac{Y^2}{(qK)^2}.
\end{equation}
\end{Le}
The proof of Lemma \ref{BK lemma 7.1} given in the article~\cite{BK} is based on the paper~\cite{BKS}. In ~\cite{BKS} using the spectral theory of automorphic forms statements similar to Lemma \ref{BK lemma 7.1} are proved. However, the main purpose of this article is to present an elementary proof of the Bourgain-Kontorovich's  result. This accounts for the fact that we prefer not to use Lemma \ref{BK lemma 7.1} in the proof of our main Theorem \ref{main}. Instead of it we use another result, although a weaker one, obtained in our paper with the use of the estimates of Kloosterman sums. We note that the problem of obtaining an elementary proof of the Bourgain-Kontorovich's theorem  was stated by Moshchevitin in ~\cite{NG}. The main result of the paper is the following theorem.
\begin{Th}\label{main}
For any alphabet $\A$ with
\begin{equation}\label{KFcondition}
\delta_{\A}>1-\frac{1}{9+\sqrt{61}}=0,9405\ldots,
\end{equation}
the following inequality (positive proportion) holds
\begin{equation}\label{KFresult}
\#\D_{\A}(N)\gg N.
\end{equation}
\end{Th}
\begin{Zam}
It is proved ~\cite{Jenkinson} that $\delta_{13}=0,9445\ldots,$ From this follows that the alphabet $\left\{1,2,\ldots,2,13\right\}$ satisfies the condition of Theorem \ref{main}. 
\end{Zam}
This article is a continuation of our paper ~\cite{FK}.  So we will heavily refer to statements and constructions in ~\cite{FK}. It should be mentioned that the proof of the Theorem \ref{main} repeats significantly  the proof of the Theorem \ref{uslov} in ~\cite{FK}. \par
In sections \ref{good function and kusic th} and \ref{method of hensley} an important accessory Theorem 6.1. from ~\cite{FK} will be proved.
\section{Обозначения}
В дальнейшем $\epsilon_0=\epsilon_0(\A)\in(0,\frac{1}{1000}).$ Для двух функций $f(x), g(x)$ используемый знак Виноградова $f(x)\ll g(x),$  означает существование константы $C,$ зависящей от $A,\epsilon_0,$ такой что $|f(x)|\le Cg(x).$  Аналогичный смысл имеет обозначение $f(x)=O(g(x)).$ Случай обоюдной оценки ${f(x)\ll g(x)\ll f(x)}$ обозначается стандартным образом $f(x)\asymp g(x).$ Используется также традиционное обозначение $e(x)=\exp(2\pi ix).$ Мощность конечного множества $S$ обозначается через $|S|$ или $\#S.$ Через $[\alpha]$ и $\|\alpha\|$ для действительных чисел $\alpha$ обозначаются, соответственно, целая часть числа $\alpha$ и расстояние от $\alpha$ до ближайшего целого:
$$[\alpha]=\max\left\{z\in\mathbb{Z}|z\le \alpha \right\},\,$$
$$\|\alpha\|=\min\left\{|z-\alpha|\Bigl|z\in\mathbb{Z}\right\}.$$
Обозначим, следуя Н.М. Коробову, через $\delta_q(a)$ характеристическую функцию делимости на натуральное число q:
\begin{gather}\label{deltasum}
\delta_q(a)=\frac{1}{q}\sum_{x=1}^q\exp\left(2\pi i\frac{ax}{q}\right)=
\left\{
              \begin{array}{ll}
                1, & \hbox{если $a\equiv 0 \pmod{q}$;} \\
                0, & \hbox{иначе.}
              \end{array}
\right.
\end{gather}
Обозначим через $\sigma_{\alpha}(q)$ следующую сумму $\sigma_{\alpha}(q)=\sum_{d|q}d^{\alpha}.$

\section{Благодарности}
Мы благодарим профессора Н.Г. Мощевитина за неоднократное обсуждение наших результатов. Именно он обратил наше внимание на статью ~\cite{BK}. Мы также выражаем благодарность И.Д. Шкредову и И.С. Резвяковой за вопросы и комментарии во время наших докладов.
\section{Оценки сумм Клоостермана и их применение.}
Для натурального $q$ и целых $l,\,m,\,n$ следуя Устинову ~\cite{Ustinov}, определим суммы Клоостермана следующим образом:
\begin{gather}\label{k1}
K_q(l,m,n)=\sum_{x,y=1}^q \delta_q(xy-l)e\left(\frac{mx+ny}{q}\right),
\end{gather}
где $\delta_q(x)$ определена в \eqref{deltasum}. При $l=1$ это определение совпадает с классическим определением сумм Клоостермана
\begin{gather}\label{k2}
K_q(m,n)=K_q(1,m,n)=\sum_{x,y=1}^q \delta_q(xy-1)e\left(\frac{mx+ny}{q}\right).
\end{gather}
Для классических сумм Клоостермана Эстерманом \cite{Estermann} было доказано неравенство
\begin{gather}\label{k3}
\left|K_q(m,n)\right|\le\sigma_0(q)(m,n,q)^{1/2}\sqrt{q}.
\end{gather}
Мы же будем использовать следующий результат, полученный Устиновым ~\cite{Ustinov}.
\begin{Le}\label{lemma-k-1}(\textbf{Устинов}~\cite{Ustinov})
Пусть $q$~--натуральное, $l,\,m,\,n$~--целые, тогда
\begin{gather}\label{k-4}
\left|K_q(l,m,n)\right|\le f_q(l,m,n)\sqrt{q},
\end{gather}
где
\begin{gather}\label{k-5}
f_q(l,m,n)=\sigma_0(q)\sigma_0((l,m,n,q))(lm,ln,mn,q)^{1/2}.
\end{gather}
\end{Le}
\begin{Zam}\label{zam-k-1}
Используя оценку $$\sigma_0(q)\ll_{\epsilon}q^{\epsilon},\quad\forall\epsilon>0,$$
из \eqref{k-5} получаем неравенство
\begin{gather}\label{k-6}
f_q(l,m,n)\ll q^{\epsilon}(l,q)^{\epsilon}(lm,ln,mn,q)^{1/2},
\end{gather}
которое мы и будем использовать.
\end{Zam}
\begin{Le}\label{lemma-k-2}(\textbf{Устинов}~\cite{Ustinov})
Пусть $q$~--натуральное, $Q_1,\,Q_2,\,P_1,\,P_2$~--вещественные, $0\le P_1,\,P_2\le q,$  тогда для суммы
\begin{gather}\label{k-7}
\Phi_q^{(0)}(Q_1,Q_2;P_1,P_2)=\sum_{Q_1<u\le Q_1+P_1\atop Q_2<v\le Q_2+P_2}\delta_q(uv\pm1)
\end{gather}
справедлива асимптотическая формула
\begin{gather}\label{k-8}
\Phi_q^{(0)}(Q_1,Q_2;P_1,P_2)=\frac{\varphi(q)}{q^2}P_1P_2+O(\psi_1(q)),
\end{gather}
где
\begin{gather}\label{k-9}
\psi_1(q)=\sigma_0(q)\log^2(q+1)q^{1/2}.
\end{gather}
\end{Le}
Следующий результат является уточнением аналогичного результата Устинова ~\cite{Ustinov}.
\begin{Le}\label{lemma-k-3}
Пусть $q$~--натуральное, $Q_1,\,Q_2$~--вещественные $0\le Q_1,\,Q_2\le q,$ $f(u)$~--невозрастающая функция на отрезке $[0;q]$ и $f(0)\le q,$  тогда
\begin{gather}\notag
\sum_{Q_1<u\le Q_2}\sum_{0<v\le f(u)}\delta_q(uv\pm1)=\frac{\varphi(q)}{q^2}\int_{Q_1}^{Q_2}f(u)du+\\+
O\left(\frac{\varphi(q)}{q^2}\frac{Q_2-Q_1}{T}(f(Q_1)-f(Q_2))\right)+O(T\psi_1(q)),\label{k-10}
\end{gather}
где $T$~--любое натуральное число.
\begin{proof}
Разобьем отрезок $[Q_1;Q_2]$ на $T$ частей
\begin{gather*}
[u_{j-1},u_{j}],\quad j=1,\ldots,T,\quad \mbox{где}\quad u_j=Q_1+j\frac{Q_2-Q_1}{T},\quad j=0,1,\ldots,T.
\end{gather*}
Обозначим
\begin{gather*}
S_j=\sum_{u_{j-1}<u\le u_{j}}\sum_{0<v\le f(u)}\delta_q(uv\pm1),\quad j=1,\ldots,T.
\end{gather*}
Заметим, что
\begin{gather}\label{k-11}
\sum_{Q_1<u\le Q_2}\sum_{0<v\le f(u)}\delta_q(uv\pm1)=\sum_{j=1}^T S_j.
\end{gather}
Ввиду не возрастания функции $f(u)$ на отрезке $[0;q],$ получаем:
\begin{gather}\label{k-11-1}
\sum_{u_{j-1}<u\le u_{j}}\sum_{0<v\le f(u_{j})}\delta_q(uv\pm1)\le
S_j\le\sum_{u_{j-1}<u\le u_{j}}\sum_{0<v\le f(u_{j-1})}\delta_q(uv\pm1).
\end{gather}
Применяя к левой и правой части неравенства \eqref{k-11-1} лемму \ref{lemma-k-2}, получаем:
\begin{gather}\label{k-11-2}
\frac{\varphi(q)}{q^2}\frac{Q_2-Q_1}{T}f(u_{j})+O(\psi_1(q))\le
S_j\le
\frac{\varphi(q)}{q^2}\frac{Q_2-Q_1}{T}f(u_{j-1})+O(\psi_1(q)).
\end{gather}
Суммируя \eqref{k-11-2} по $j=1,\ldots,T,$  получаем:
\begin{gather}\label{k-11-3}
\frac{\varphi(q)}{q^2}\sum_{j=1}^{T}\frac{Q_2-Q_1}{T}f(u_{j})+TO(\psi_1(q))\le
\sum_{j=1}^{T}S_j\le
\frac{\varphi(q)}{q^2}\sum_{j=1}^{T}\frac{Q_2-Q_1}{T}f(u_{j-1})+TO(\psi_1(q)).
\end{gather}
Ввиду не возрастания функции $f(u)$ на отрезке $[0;q],$ получаем:
\begin{gather}\label{k-11-4}
\sum_{j=1}^{T}\frac{Q_2-Q_1}{T}f(u_{j})\le\int_{Q_1}^{Q_2}f(u)du\le\sum_{j=1}^{T}\frac{Q_2-Q_1}{T}f(u_{j-1})
\end{gather}
и
\begin{gather}\label{k-11-5}
\sum_{j=1}^{T}\frac{Q_2-Q_1}{T}f(u_{j-1})-\sum_{j=1}^{T}\frac{Q_2-Q_1}{T}f(u_{j})=
\frac{Q_2-Q_1}{T}\left(f(Q_1)-f(Q_2)\right).
\end{gather}
Следовательно, из  \eqref{k-11-3} и соотношений  \eqref{k-11-4}, \eqref{k-11-5} получаем:
\begin{gather}\label{k-11-6}
\sum_{j=1}^{T}S_j=
\frac{\varphi(q)}{q^2}\int_{Q_1}^{Q_2}f(u)du+
O\left(\frac{\varphi(q)}{q^2}\frac{Q_2-Q_1}{T}(f(Q_1)-f(Q_2))\right)+O(T\psi_1(q)).
\end{gather}
Подставляя \eqref{k-11-6} в \eqref{k-11}, получаем \eqref{k-10}. Лемма доказана.
\end{proof}
\end{Le}
Докажем две леммы, обобщающие лемму \ref{lemma-k-2}. Их доказательство идейно повторяет доказательство леммы \ref{lemma-k-2} ~\cite{Ustinov}. В дальнейшем  запись вида $\mathop{{\sum}'}_{-s/2<m\le s/2}$ означает, что суммирование ведется по $m\neq0.$
\begin{Le}\label{lemma-k-4}
Пусть $q,\,d,\,k$~--натуральные, $s=qd,$ $Q_1,\,Q_2,\,P_1,\,P_2$~--вещественные, $0\le P_1\le d,\,0\le P_2\le s,$  тогда для суммы
\begin{gather}\label{k-12}
\Phi_{s,q,k}^{(1)}(Q_1,Q_2;P_1,P_2)=\sum_{Q_1<u\le Q_1+P_1\atop Q_2<v\le Q_2+P_2}\delta_s(uv\pm q)\delta_q(v)
e\left(\frac{k}{q}u\right)
\end{gather}
справедлива асимптотическая формула
\begin{gather}\label{k-13}
\Phi_{s,q,k}^{(1)}(Q_1,Q_2;P_1,P_2)=\frac{P_1P_2}{qs^2}\sum_{l=1}^q K_s(\pm q,dl,dk)+
O(\sum_{j=0}^3R_j(k)),
\end{gather}
где
\begin{gather}\label{k-14}
R_0(k)=\frac{1}{qs}\left|\sum_{l=1}^q K_s(\pm q,dl,dk)\right|,
\end{gather}
\begin{gather}\label{k-15}
R_1(k)=\frac{1}{q^2}\mathop{{\sum}'}_{-s/2<m\le s/2}\frac{1}{|m|}\left|\sum_{l=1}^qK_s(\pm q,m+dl,dk)\right|,
\end{gather}
\begin{gather}\label{k-16}
R_2(k)=\frac{1}{q}\mathop{{\sum}'}_{-s/2<n\le s/2}\frac{1}{|n|}\left|\sum_{l=1}^qK_s(\pm q,dl,n+dk)\right|,
\end{gather}
\begin{gather}\label{k-17}
R_3(k)=\frac{1}{q}\mathop{{\sum}'}_{-s/2<m\le s/2}\mathop{{\sum}'}_{-s/2<n\le s/2}\frac{1}{|mn|}\left|\sum_{l=1}^qK_s(\pm q,m+dl,n+dk).
\right|
\end{gather}
\begin{proof}
Положим
\begin{gather*}
M_1=[Q_1],\quad N_1=[Q_1+P_1]-[Q_1],
M_2=[Q_2],\quad N_2=[Q_2+P_2]-[Q_2],
\end{gather*}
тогда
\begin{gather*}
\Phi_{s,q,k}^{(1)}(Q_1,Q_2;P_1,P_2)=\Phi_{s,q,k}^{(1)}(M_1,M_2;N_1,N_2).
\end{gather*}
Очевидно, что
\begin{gather}\label{k-18}
N_j=P_j-\left\{Q_j+P_j\right\}+\left\{Q_j\right\}, \quad j=1,2
\end{gather}
и, следовательно $0\le N_1\le d,\,0\le N_2\le s.$  Определим две характеристические функции  для $M_j<x\le M_j+s,\,j=1,2:$
\begin{gather}
\chi_j(x)=
\left\{
              \begin{array}{ll}
                1, & \hbox{если $M_j<x\le M_j+N_j$;} \\
                0, & \hbox{если $M_j+N_j<x\le M_j+s$.}
              \end{array}
\right.\label{k-18-1}
\end{gather}
Тогда их разложение в конечный ряд Фурье имеет следующий вид:
\begin{gather}\label{k-19}
\chi_j(x)=\sum_{-s/2<n\le s/2}\hat{\chi}_j(n)e\left(\frac{nx}{s}\right),\quad\mbox{где}\quad
\hat{\chi}_j(n)=\frac{1}{s}\sum_{x=M_j+1}^{M_j+N_j}e\left(\frac{-nx}{s}\right).
\end{gather}
При $-\frac{s}{2}<n\le\frac{s}{2},\,n\neq0,$ после суммирования геометрической прогрессии получаем:
\begin{gather}\label{k-19-1}
\left|\hat{\chi}_j(n)\right|=\frac{1}{s}\frac{|1-e(\frac{nN_j}{s})|}{|1-e(\frac{n}{s})|}\le
\frac{1}{s}\frac{1}{|\sin(\pi n/s)|}\le\frac{1}{2|n|}.
\end{gather}
Следовательно,
\begin{gather}\label{k-20}
\hat{\chi}_j(0)=\frac{N_j}{s},\quad
\left|\hat{\chi}_j(n)\right|\le\frac{1}{2|n|}\quad\mbox{при}\quad -\frac{s}{2}<n\le\frac{s}{2},\, n\neq0.
\end{gather}
Из определения \eqref{k-18-1} функций $\chi_j(x), j=1,2,$ следует, что
\begin{gather}\label{k-20-1}
\Phi_{s,q,k}^{(1)}(M_1,M_2;N_1,N_2)=
\sum_{u=M_1+1}^{M_1+s}\sum_{v=M_2+1}^{M_2+s}\delta_s(uv\pm q)\delta_q(v)\chi_1(u)\chi_2(v)
e\left(\frac{k}{q}u\right).
\end{gather}
Подставляя в \eqref{k-20-1} разложения в ряд Фурье функций $\chi_1(u),\,\chi_2(v),$ (см. \eqref{k-19}), получаем:
\begin{gather*}
\Phi_{s,q,k}^{(1)}(M_1,M_2;N_1,N_2)=
\sum_{-s/2<m,n\le s/2}\hat{\chi}_1(n)\hat{\chi}_2(m)
\sum_{u=1}^{s}\sum_{v=1}^{s}\delta_s(uv\pm q)\delta_q(v)e\left(\frac{nu+mv}{s}\right)
e\left(\frac{k}{q}u\right).
\end{gather*}
Используя определение \eqref{deltasum}, получаем:
\begin{gather*}
\Phi_{s,q,k}^{(1)}(M_1,M_2;N_1,N_2)=\frac{1}{q}
\sum_{-s/2<m,n\le s/2}\hat{\chi}_1(n)\hat{\chi}_2(m)\sum_{l=1}^{q}
\sum_{u=1}^{s}\sum_{v=1}^{s}\delta_s(uv\pm q)e\left(\frac{nu+mv}{s}+\frac{ku}{q}+\frac{lv}{q}\right).
\end{gather*}
По условию $s=qd,$ следовательно,
\begin{gather}\label{k-20-1-1}
\frac{nu+mv}{s}+\frac{ku}{q}+\frac{lv}{q}=\frac{(n+dk)u+(m+dl)v}{s}.
\end{gather}
Используя \eqref{k-20-1-1} и определение \eqref{k1}, получаем:
\begin{gather*}
\Phi_{s,q,k}^{(1)}(M_1,M_2;N_1,N_2)=\frac{1}{q}
\sum_{-s/2<m,n\le s/2}\hat{\chi}_1(n)\hat{\chi}_2(m)\sum_{l=1}^{q}
K_s(\pm q,m+dl,n+dk).
\end{gather*}
Выделяя слагаемое с $m=n=0$ и слагаемые с $m=0$ или $n=0,$ используя \eqref{k-20}, получаем:
\begin{gather*}
\Phi_{s,q,k}^{(1)}(M_1,M_2;N_1,N_2)=
\frac{N_1N_2}{qs^2}\sum_{l=1}^q K_s(\pm q,dl,dk)+
\frac{N_1}{qs}\mathop{{\sum}'}_{-s/2<m\le s/2}\hat{\chi}_2(m)\sum_{l=1}^qK_s(\pm q,m+dl,dk)+\\+
\frac{N_2}{qs}\mathop{{\sum}'}_{-s/2<n\le s/2}\hat{\chi}_1(n)\sum_{l=1}^qK_s(\pm q,dl,n+dk)+\\+
\frac{1}{q}\mathop{{\sum}'}_{-s/2<m\le s/2}\mathop{{\sum}'}_{-s/2<n\le s/2}\hat{\chi}_1(n)\hat{\chi}_2(m)\sum_{l=1}^qK_s(\pm q,m+dl,n+dk).
\end{gather*}
Из \eqref{k-18} следует, что $P_j-1\le N_j\le P_j+1$ для $j=1,2.$ Используя это и оценки коэффициентов Фурье из  \eqref{k-20}, получаем
\begin{gather*}
\Phi_{s,q,k}^{(1)}(M_1,M_2;N_1,N_2)
=\frac{N_1N_2}{qs^2}\sum_{l=1}^q K_s(\pm q,dl,dk)+O(R_1(k))+O(R_2(k))+O(R_3(k))
\end{gather*}
Далее, так как
\begin{gather*}
\left|P_1P_2-N_1N_2\right|=\left|P_1(P_2-N_2)+N_2(P_1-N_1)\right|\le P_1+N_2\le P_1+P_2+1\ll s,
\end{gather*}
то получаем:
\begin{gather*}
\Phi_{s,q,k}^{(1)}(M_1,M_2;N_1,N_2)
=\frac{P_1P_2}{qs^2}\sum_{l=1}^q K_s(\pm q,dl,dk)+O(\sum_{j=0}^3R_j(k)).
\end{gather*}
Лемма доказана.
\end{proof}
\end{Le}

Следующая лемма является обобщением леммы \ref{lemma-k-3}.
\begin{Le}\label{lemma-k-6}
Пусть $q,\,d,\,k$~--натуральные, $s=qd.$ Пусть $Q_1,\,Q_2$~--вещественные, ${0\le Q_1,\,Q_2\le d,}$ $f(u)$~--невозрастающая функция на отрезке $[0;s]$ и $f(0)\le s,$  тогда
\begin{gather}\notag
\sum_{Q_1<u\le Q_2}\sum_{0<v\le f(u)}\delta_s(uv\pm q)\delta_q(v)e\left(\frac{k}{q}u\right)=
\frac{1}{qs^2}\sum_{l=1}^q K_s(\pm q,dl,dk)\int_{Q_1}^{Q_2}f(u)du+\\+
O\left(\frac{1}{qs^2}\left|\sum_{l=1}^q K_s(\pm q,dl,dk)\right|\frac{Q_2-Q_1}{T}(f(Q_1)-f(Q_2))\right)+
O(T\sum_{j=0}^3R_j(k)),\label{k-27}
\end{gather}
где $T$~--любое натуральное число,а величины $R_0(k),R_1(k),R_2(k),R_3(k)$ определены в \eqref{k-14}~--\eqref{k-17}, соответственно.
\begin{proof}
Доказательство этой леммы, опирающееся на лемму \ref{lemma-k-4}, дословно повторяет вывод леммы \ref{lemma-k-3} из леммы \ref{lemma-k-2}.
\end{proof}
\end{Le}
\begin{Zam}\label{zam-k-2}
Утверждение, аналогичное лемме \ref{lemma-k-4}, может быть доказано для сумм следующего вида
\begin{gather}\label{k-21}
\Phi_{s,q,k}^{(2)}(Q_1,Q_2;P_1,P_2)=\sum_{Q_1<u\le Q_1+P_1\atop Q_2<v\le Q_2+P_2}\delta_s(uv\pm q^2)\delta_q(u)\delta_q(v)
e\left(\frac{p}{q^2}u+\frac{t}{q^2}v\right),
\end{gather}
где $s=q^2d.$ На основе \eqref{k-21} может быть получена формула, аналогичная \eqref{k-27}, для следующей суммы:
\begin{gather}\label{k-28}
\sum_{Q_1<u\le Q_2}\sum_{0<v\le f(u)}
\delta_s(uv\pm q^2)\delta_q(u)\delta_q(v)e\left(\frac{p}{q^2}u+\frac{t}{q^2}v\right).
\end{gather}
\end{Zam}

\section{Видоизменение леммы \ref{BK lemma 7.1}.}
Пусть даны $X>Y>0,\,K>0$~--вещественные, $q\asymp Q$~--натуральное. Пусть так же даны два вектора
$$\eta=(x,y)^t,\,\eta'=(u,v)^t$$
с натуральными координатами, удовлетворяющими соотношениям
\begin{gather}\label{f1}
\frac{y}{A}\le x\le y,\quad \frac{X}{Y}\ll y\ll\frac{X}{Y},\quad (x,y)=1
\end{gather}
\begin{gather}\label{f2}
\frac{v}{A}\le u\le v,\quad X\ll v\ll X,\,\frac{2X}{K}\le v,\quad (u,v)=1,
\end{gather}
причем все константы в знаках $\ll$~--абсолютные. Рассмотрим множество матриц
\begin{gather}\label{f3}
M=\left\{
\gamma=
\begin{pmatrix}
a & b \\
c & d
\end{pmatrix}\Biggl|\quad
0\le a\le b,\,1\le c\le d,\, 1\le b\le d,\,det\gamma=1
\right\}.
\end{gather}
В данном параграфе речь пойдет об оценке мощности множества матриц
\begin{gather}\label{f4}
\M=\left\{
\gamma\in M\Biggl|\quad
\left|\gamma\eta-\eta'\right|\le\frac{X}{K},\,\gamma\eta\equiv\eta'\pmod{q},\, \frac{Y}{\kappa_1}\le\|\gamma\|\le Y
\right\},
\end{gather}
где $\kappa_1$~--какая-либо константа. Рассмотрим подробнее условия на $\gamma\in\M.$ Условие $\left|\gamma\eta-\eta'\right|\le\frac{X}{K}$ равносильно:
\begin{gather*}
-\frac{X}{K}<ax+by-u<\frac{X}{K},\quad -\frac{X}{K}<cx+dy-v<\frac{X}{K}.
\end{gather*}
По определению, $\|\gamma\|=d,$ следовательно, мощность множества $\M$ не превосходит числа решений $(a,b,c,d)$ системы
\begin{gather}
\left\{
  \begin{array}{ll}
    -\frac{X}{K}<ax+by-u<\frac{X}{K}, &  \\
    -\frac{X}{K}<cx+dy-v<\frac{X}{K}, & \\
    ax+by\equiv u\pmod{q}, &  \\
    cx+dy\equiv v\pmod{q}, &  \\
    \frac{Y}{\kappa_1}\le d\le Y, &  \\
    \gamma\in M.
    \end{array}
\right.\label{f5}
\end{gather}
Так как $\gamma\in M,$ то $det\gamma=ad-bc=1$ и, следовательно,
$$a=\frac{bc+1}{d},\quad d|(bc+1).$$
Подставляя $a$ в первое неравенство системы \eqref{f5}, получаем:
\begin{gather*}
 -d\frac{X}{K}<b(cx+dy)+ x-du<d\frac{X}{K}.
\end{gather*}
Далее,
\begin{gather*}
d\left(u-\frac{X}{K}\right)-x<b(cx+dy)<d\left(u+\frac{X}{K}\right)-x.
\end{gather*}
Из определения множества $M$ получаем, что $c>0,\,d>0;$ ввиду того что у вектора $\eta=(x,y)^t$ натуральные координаты, получаем $0<cx+dy,$ следовательно
\begin{gather}\label{f6}
\frac{d\left(u-\frac{X}{K}\right)-x}{cx+dy}<b<\frac{d\left(u+\frac{X}{K}\right)-x}{cx+dy}.
\end{gather}
Из второго неравенства системы \eqref{f5} получаем:
\begin{gather}\label{f7}
\frac{\left(v-\frac{X}{K}\right)-dy}{x}<c<\frac{d\left(v+\frac{X}{K}\right)-dy}{x}.
\end{gather}
Обозначим:
\begin{gather}\label{f8}
c_1=\frac{\left(v-\frac{X}{K}\right)-dy}{x},\quad
c_2=\frac{d\left(v+\frac{X}{K}\right)-dy}{x},
\end{gather}
\begin{gather}\label{f9}
f_1(c)=\frac{d\left(u-\frac{X}{K}\right)-x}{cx+dy},\quad
f_2(c)=\frac{d\left(u+\frac{X}{K}\right)-x}{cx+dy}.
\end{gather}
Теперь система \eqref{f5} принимает вид:
\begin{gather}
\left\{
  \begin{array}{ll}
    c_1<c<c_2,\,1\le c\le d,&  \\
    f_1(c)<b<f_2(c),\, 1\le b\le d & \\
    ax+by\equiv u\pmod{q}, &  \\
    cx+dy\equiv v\pmod{q}, &  \\
    a=\frac{bc+1}{d},\quad d|(bc+1), &  \\
    \frac{Y}{\kappa_1}\le d\le Y.
    \end{array}
\right.\label{f10}
\end{gather}
Две первые строки этой системы задают область
\begin{gather}\label{f11}
\Theta=\left\{(c,b)\Biggl| \tilde{c}_1<c<\tilde{c}_2,\, \tilde{f}_1(c)<b<\tilde{f}_2(c)\right\},
\end{gather}
причем функции $\tilde{f}_1(c),\,\tilde{f}_2(c)$ не возрастают на отрезке $[\tilde{c}_1;\tilde{c}_2].$ В дальнейшем мы часто будем использовать следующие обозначения:
\begin{gather}\label{f12}
I_{c}=[ \tilde{c}_1,\tilde{c}_2],\, I_{b}=[\tilde{f}_1(c),\tilde{f}_2(c)].
\end{gather}
Мы оценим $\left|\Theta\right|$ следующим образом:
\begin{gather*}
\left|\Theta\right|\le\int_{c_1}^{c_2}\left(f_2(s)-f_1(s)\right)ds.
\end{gather*}
Подставляя $c_1,c_2,f_1(c),f_2(c)$ из \eqref{f8} и \eqref{f9}, производя замену $z=sx+dy,$ получаем:
\begin{gather}\label{f13}
\left|\Theta\right|\le\int_{c_1}^{c_2}\frac{2d\frac{X}{K}}{sx+dy}ds=
2d\frac{X}{K}\int_{v-\frac{X}{K}}^{v+\frac{X}{K}}\frac{1}{xz}dz=
2d\frac{X}{K}\frac{1}{x}\log\frac{1+\frac{X}{vK}}{1-\frac{X}{vK}}.
\end{gather}
Ввиду того, что
\begin{gather*}
\log\frac{1+z}{1-z}\ll z,\; \mbox{при}\; 0<z<1,
\end{gather*}
а  $0<\frac{X}{vK}\le\frac{1}{2},$ получаем:
\begin{gather*}
\left|\Theta\right|\ll d\frac{X^2}{K^2}\frac{1}{vx}.
\end{gather*}
Ввиду \eqref{f1} и \eqref{f2} получаем:
\begin{gather}\label{f14}
\left|\Theta\right|\ll d\frac{Y}{K^2}.
\end{gather}
Заметим, что
\begin{gather}\label{f15}
\tilde{c}_2-\tilde{c}_1\le c_2-c_1=\frac{X}{K}\frac{2}{x}\ll\frac{Y}{K},
\end{gather}
где в последнем переходе мы воспользовались \eqref{f1}. Нам так же необходимо оценить разности $\tilde{f}_{i}(\tilde{c}_1)-\tilde{f}_{i}(\tilde{c}_2).$
Оценим $\tilde{f}_{i}(\tilde{c}_1)-\tilde{f}_{i}(\tilde{c}_2)$ следующим образом:
\begin{gather}\label{f16}
\tilde{f}_{i}(\tilde{c}_1)-\tilde{f}_{i}(\tilde{c}_2)\le f_{i}(c_1)-f_{i}(c_2),\quad i=1,2.
\end{gather}
Ввиду \eqref{f8} и \eqref{f9} получаем:
\begin{gather*}
f_1(c_1)=\frac{d\left(u-\frac{X}{K}\right)-x}{v-\frac{X}{K}},\quad
f_1(c_2)=\frac{d\left(u-\frac{X}{K}\right)-x}{v+\frac{X}{K}}.
\end{gather*}
Следовательно,
\begin{gather*}
f_1(c_1)-f_1(c_2)=\frac{d\left(u-\frac{X}{K}\right)-x}{v^2-\frac{X^2}{K^2}}\frac{2X}{K}\le\frac{2dX}{K}
\frac{v}{v^2-\frac{X^2}{K^2}}\le\frac{2dX}{K}\frac{1}{v-\frac{X}{K}}.
\end{gather*}
Используя \eqref{f1} и \eqref{f2}, получаем:
\begin{gather}\label{f17}
f_1(c_1)-f_{1}(c_2)\ll\frac{d}{K}.
\end{gather}
Совершенно аналогично получаем:
\begin{gather}\label{f18}
f_2(c_1)-f_{2}(c_2)\ll\frac{d}{K}.
\end{gather}
Подставляя \eqref{f17} и \eqref{f18} в \eqref{f16}, получаем:
\begin{gather}\label{f19}
\tilde{f}_{i}(\tilde{c}_1)-\tilde{f}_{i}(\tilde{c}_2)\ll\frac{d}{K},\quad i=1,2.
\end{gather}
Теперь все готово для того, чтобы оценить мощность множества $\M.$



\begin{Th}\label{theoremF-2}
Если $X>Y>K^{4}q^3,$ то для любого $\epsilon>0$ справедлива следующая оценка мощности множества $\M:$
\begin{gather}\label{f-2-1}
\left|\M\right|\ll_{\epsilon}\frac{Y^{2+\epsilon}}{K^2}\frac{q^{\epsilon}}{\sqrt{q}}.
\end{gather}
\begin{proof}
Очевидно, что число решений системы \eqref{f10} не превосходит числа решений $\sys=S(u,v,x,y,q)$ следующей системы:
\begin{gather}
\left\{
  \begin{array}{ll}
    c_1<c<c_2,\,1\le c\le d,&  \\
    f_1(c)<b<f_2(c),\, 1\le b\le d & \\
    cx+dy\equiv v\pmod{q}, &  \\
    d|(bc+1), &  \\
    \frac{Y}{\kappa_1}\le d\le Y,
    \end{array}
\right.\label{f-2-2}
\end{gather}
в переменных $(b,c,d).$ Чтобы не вводить дополнительных обозначений для остатков от деления $x,y,v$ на $q,$ будем считать, что $0<x,y,v\le q.$ Обозначим
\begin{gather}\label{f-2-2.1}
\alpha_1=(q,x),\quad \alpha_2=(q,v), \alpha_3=\frac{\alpha_1\alpha_2}{(\alpha_1,\alpha_2)}.
\end{gather}
Так как $dy-v\equiv -cx\pmod{q},$ то $\alpha_1|(dy-v).$ Ввиду $(x,y)=1$ получаем $(\alpha_1,y)=1.$ Тогда из
$dy\equiv v\pmod{\alpha_1}$ следует, что $d\equiv v'\pmod{\alpha_1},$ где $0<v'\le\alpha_1.$ Поэтому число решений системы \eqref{f-2-2} равно
\begin{gather}\label{f-2-3}
\sys=\sum_{\frac{Y}{\kappa_1}\le d\le Y}\delta_{\alpha_1}(d-v')\sum_{c\in I_c}\sum_{b\in I_b}\delta_d(bc+1)\delta_q(cx+dy-v).
\end{gather}
Так как
\begin{gather*}
\delta_q(cx+dy-v)=\frac{1}{q}\sum_{k=1}^q e\left(\frac{cx+dy-v}{q}k\right),
\end{gather*}
то из \eqref{f-2-3} получаем, что
\begin{gather}\label{f-2-4}
\sys\le
\frac{1}{q}\sum_{k=1}^q
\sum_{\frac{Y}{\kappa_1}\le d\le Y}\delta_{\alpha_1}(d-v')
\left|\sum_{c\in I_c}\sum_{b\in I_b}\delta_d(bc+1)e\left(\frac{kx}{q}c\right)\right|.
\end{gather}
Если $bc\equiv-1\pmod{d},$ то $(qb)c\equiv-q\pmod{qd}.$ Положив $s=qd,$ получаем:
\begin{gather}\label{f-2-5}
\sum_{c\in I_c}\sum_{b\in I_b}\delta_d(bc+1)e\left(\frac{kx}{q}c\right)=
\sum_{c\in I_c}\sum_{b\in qI_b}\delta_s(bc+q)\delta_q(b)e\left(\frac{kx}{q}c\right).
\end{gather}
Учитывая, что $I_c=[\tilde{c}_1,\tilde{c}_2],\, qI_b=[q\tilde{f}_1(c),q\tilde{f}_2(c)],$ применяем лемму \ref{lemma-k-6} с $k$ равным $kx,$ получаем:
\begin{gather}\notag
\sum_{c\in I_c}\sum_{b\in I_b}\delta_d(bc+1)e\left(\frac{kx}{q}c\right)=
\sum_{c\in I_c}\sum_{b\in qI_b}\delta_s(bc+q)\delta_q(b)e\left(\frac{kx}{q}c\right)=\\\notag=
\frac{1}{qs^2}\sum_{l=1}^q K_s(q,dl,dkx)
\int_{\tilde{c}_1}^{\tilde{c}_2}(q\tilde{f}_2(s)-q\tilde{f}_1(s))ds+\\\notag+
O\left(\frac{1}{qs^2}\left|\sum_{l=1}^q K_s(q,dl,dkx)\right|
\frac{\tilde{c}_2-\tilde{c}_1}{T}(q\tilde{f}_1(c)-q\tilde{f}_2(c))\right)+O(T\sum_{j=0}^3R_j(kx)).\label{f-2-6}
\end{gather}
Из определения $\Theta$ и оценки \eqref{f14} получаем:
\begin{gather}\label{f-2-6-1}
\frac{1}{q}\int_{\tilde{c}_1}^{\tilde{c}_2}(q\tilde{f}_2(s)-q\tilde{f}_1(s))ds=|\Theta|\ll d\frac{Y}{K^2}.
\end{gather}
Подставляя в \eqref{f-2-6} оценки из \eqref{f-2-6-1}, \eqref{f15} и \eqref{f19}, получаем:
\begin{gather}\notag
\left|\sum_{c\in I_c}\sum_{b\in I_b}\delta_d(bc\pm1)e\left(\frac{kx}{q}c\right)\right|\ll
\frac{d}{s^2}\frac{Y}{K^2}\sum_{l=1}^q \left|K_s(\pm q,dl,dkx)\right|+\\\notag+
O\left(\frac{1}{s^2}\left|\sum_{l=1}^q K_s(\pm q,dl,dkx)\right|
\frac{dY}{TK^2}\right)+O(T\sum_{j=0}^3R_j(kx)).\label{f-2-7}
\end{gather}
Подставляя \eqref{f-2-7} в \eqref{f-2-4}, получаем, что число решений системы \eqref{f-2-2} не превосходит
\begin{gather*}
\sys\ll\sum_{\frac{Y}{\kappa_1}\le d\le Y}\delta_{\alpha_1}(d-v')\frac{d}{qs^2}\frac{Y}{K^2}\sum_{k=1}^q
\sum_{l=1}^q \left|K_s(\pm q,dl,dkx)\right|+\\+
O\left(\sum_{\frac{Y}{\kappa_1}\le d\le Y}
\delta_{\alpha_1}(d-v')\frac{d}{qs^2}\frac{Y}{TK^2}\sum_{k=1}^q\left|\sum_{l=1}^q K_s(\pm q,dl,dkx)\right|
\right)+
\frac{T}{q}\sum_{j=0}^3\sum_{k=1}^q\sum_{\frac{Y}{\kappa_1}\le d\le Y}
\delta_{\alpha_1}(d-v')O(R_j(kx)).\label{f-2-8}
\end{gather*}
Обозначим
\begin{gather}\label{f-2-9}
\Sigma_1=\sum_{\frac{Y}{\kappa_1}\le d\le Y}\delta_{\alpha_1}(d-v')\frac{d}{qs^2}\frac{Y}{K^2}\sum_{k=1}^q
\sum_{l=1}^q \left|K_s(\pm q,dl,dkx)\right|,
\end{gather}
\begin{gather}\label{f-2-10}
\Sigma_2=O\left(\sum_{\frac{Y}{\kappa_1}\le d\le Y}
\delta_{\alpha_1}(d-v')\frac{d}{qs^2}\frac{Y}{TK^2}\sum_{k=1}^q\left|\sum_{l=1}^q K_s(\pm q,dl,dkx)\right|,
\right)
\end{gather}
\begin{gather}\label{f-2-11}
\Sigma_j=\frac{T}{q}\sum_{k=1}^q\sum_{\frac{Y}{\kappa_1}\le d\le Y}\delta_{\alpha_1}(d-v')O(R_{j-3}(kx)),\quad j=3,4,5,6.
\end{gather}
Следовательно,
\begin{gather*}
\sys\le\Sigma_1+\Sigma_2+\sum_{j=3}^6\Sigma_{j},
\end{gather*}
Подставляя оценки из формул \eqref{f-2-30},\,\eqref{f-2-31},\,\eqref{f-2-33},\,\eqref{f-2-40},\eqref{f-2-42},\,\eqref{f-2-47}, получаем что число решений системы \eqref{f-2-2} не превосходит
\begin{gather}\label{f-2-48}
\sys\ll_{\epsilon} \frac{Y^{2+\epsilon}}{K^2}q^{-1/2+\epsilon}+\frac{1}{T}\frac{Y^{2+\epsilon}}{K^2}q^{-1/2+\epsilon}+
Tq^{1+\epsilon}Y^{3/2+\epsilon}.
\end{gather}
Полагая
\begin{gather}\label{f-2-48}
T=\left[\frac{Y^{1/4}}{Kq^{3/4}}\right]+1,
\end{gather}
и учитывая, что $Y>K^{4}q^3,$ получаем формулу \eqref{f-2-1}. Для того, чтобы завершить доказательство теоремы нам необходимо оценить каждую величину $\Sigma_j,\,j=1,\ldots,6.$
\begin{enumerate}
  \item Оценим $\Sigma_1.$
\begin{Le}\label{lemma-f-2-1}
Справедлива следующая оценка
\begin{gather}\label{f-2-12}
\sum_{k=1}^q
\sum_{l=1}^q \left|K_s(\pm q,dl,dkx)\right|\ll_{\epsilon}(q,dx)^{1/2}q^{5/2+\epsilon}d^{1+\epsilon}
\end{gather}
для $\forall\epsilon>0.$
\begin{proof}
Применяя к $\left|K_s(\pm q,dl,dkx)\right|$ лемму \ref{lemma-k-1} и  используя замечание \ref{zam-k-1}, получим, что для $\forall\epsilon>0$ выполнено
\begin{gather}\label{f-2-13}
\left|K_s(\pm q,dl,dkx)\right|\le
s^{\epsilon}(q,s)^{\epsilon}(qdl,qdkx,d^2lkx,s)^{1/2}\sqrt{s}.
\end{gather}
Ввиду того, что $s=qd,$ имеем
\begin{gather*}
(qdl,qdkx,d^2lkx,s)=(sl,skx,d^2lkx,s)=(d^2lkx,qd)=d(dlkx,q).
\end{gather*}
Тогда из \eqref{f-2-13} следует, что
\begin{gather}\label{f-2-14}
\left|K_s(\pm q,dl,dkx)\right|\le
d^{\epsilon}q^{\epsilon}(dlkx,q)^{1/2}\sqrt{sd}.
\end{gather}
Подставляя \eqref{f-2-14} в \eqref{f-2-12}, получим
\begin{gather}\label{f-2-15}
\sum_{k=1}^q
\sum_{l=1}^q \left|K_s(\pm q,dl,dkx)\right|\le d^{\epsilon}q^{\epsilon}d\sqrt{q}
\sum_{k=1}^q\sum_{l=1}^q (dlkx,q)^{1/2}.
\end{gather}
Обозначим $\alpha=(q,dx)$ и $q_1=\frac{q}{\alpha}.$ Тогда $(dlkx,q)=\alpha(kl,q_1)$ и, следовательно,
\begin{gather}\label{f-2-16}
\sum_{k=1}^q\sum_{l=1}^q (dlkx,q)^{1/2}=\alpha^{1/2}\sum_{k=1}^q\sum_{l=1}^q(kl,q_1)^{1/2}.
\end{gather}
Оценим сумму из правой части \eqref{f-2-16}. Для $\forall\epsilon>0$ выполнено
\begin{gather}\notag
\sum_{k=1}^q\sum_{l=1}^q(kl,q_1)^{1/2}\le
\sum_{\gamma|q_1}\gamma^{1/2}\sum_{k=1}^q\sum_{l=1}^q\delta_{\gamma}(kl)=\\=
\sum_{\gamma|q_1}\gamma^{1/2}\sum_{\gamma_1|\gamma}\sum_{k=1}^q\delta_{\gamma_1}(k)\sum_{l=1}^q\delta_{\gamma/\gamma_1}(l)\le
\sum_{\gamma|q_1}\gamma^{1/2}\sum_{\gamma_1|\gamma}\frac{q^2}{\gamma}=q^2\sum_{\gamma|q_1}\frac{\sigma_0(\gamma)}{\gamma^{1/2}}\ll_{\epsilon}
q^{2+\epsilon}.\label{f-2-17}
\end{gather}
Подставляя \eqref{f-2-17} в \eqref{f-2-16}, получаем
\begin{gather}\label{f-2-18}
\sum_{k=1}^q\sum_{l=1}^q (dlkx,q)^{1/2}=(q,dx)^{1/2}q^{2+\epsilon}.
\end{gather}
Подставляя \eqref{f-2-18} в \eqref{f-2-15}, получаем утверждение леммы.
\end{proof}
\end{Le}
Используя лемму \ref{lemma-f-2-1} докажем следующее утверждение.
\begin{Le}\label{lemma-f-2-1-1}
Справедлива следующая оценка
\begin{gather}\label{f-2-30}
\Sigma_1\ll_{\epsilon}\frac{Y^{2+\epsilon}}{K^2}q^{-1/2+\epsilon}
\end{gather}
для $\forall\epsilon>0.$
\begin{proof}
Подставляя \eqref{f-2-12} в \eqref{f-2-9} и учитывая, что $s=qd,$ получаем:
\begin{gather}\label{f-2-19}
\Sigma_1\ll_{\epsilon}\frac{Y}{K^2}q^{-1/2+\epsilon}
\sum_{\frac{Y}{\kappa_1}\le d\le Y}\delta_{\alpha_1}(d-v')(q,dx)^{1/2}d^{\epsilon}.
\end{gather}
Из определения \eqref{f-2-2.1} получаем, что $(q,dx)=\alpha_1(\frac{q}{\alpha_1},d),$ следовательно
\begin{gather}\label{f-2-20}
\Sigma_1\ll_{\epsilon}\frac{Y^{1+\epsilon}}{K^2}q^{-1/2+\epsilon}\alpha_1^{1/2}
\sum_{\frac{Y}{\kappa_1}\le d\le Y}\delta_{\alpha_1}(d-v')\left(\frac{q}{\alpha_1},d\right)^{1/2}.
\end{gather}
Преобразуем сумму из правой части  \eqref{f-2-20}
\begin{gather}\label{f-2-21}
\sum_{\frac{Y}{\kappa_1}\le d\le Y}\delta_{\alpha_1}(d-v')\left(\frac{q}{\alpha_1},d\right)^{1/2}\le
\sum_{\gamma|\frac{q}{\alpha_1}}\gamma^{1/2}\sum_{\frac{Y}{\kappa_1}\le d\le Y}\delta_{\alpha_1}(d-v')\delta_{\gamma}(d).
\end{gather}
Оценим сумму из правой части  \eqref{f-2-21}.
\begin{gather}\notag
\sum_{\frac{Y}{\kappa_1}\le d\le Y}\delta_{\alpha_1}(d-v')\delta_{\gamma}(d)=
\sum_{\frac{Y}{\gamma\kappa_1}\le d\le \frac{Y}{\gamma}}\delta_{\alpha_1}(d\gamma-v')=
\frac{1}{\alpha_1}\sum_{\frac{Y}{\gamma\kappa_1}\le d\le \frac{Y}{\gamma}}
\sum_{n=1}^{\alpha_1}e\left(\frac{d\gamma-v'}{\alpha_1}n\right)
=\\=\frac{1}{\alpha_1}\sum_{n=1}^{\alpha_1}e\left(\frac{-v'}{\alpha_1}n\right)
\sum_{\frac{Y}{\gamma\kappa_1}\le d\le \frac{Y}{\gamma}}e\left(\frac{n\gamma}{\alpha_1}d\right).\label{f-2-22}
\end{gather}
Оценивая правую часть \eqref{f-2-22} по модулю, имеем
\begin{gather}\label{f-2-23}
\sum_{\frac{Y}{\kappa_1}\le d\le Y}\delta_{\alpha_1}(d-v')\delta_{\gamma}(d)\le
\frac{1}{\alpha_1}\sum_{n=1}^{\alpha_1}\left|\sum_{\frac{Y}{\gamma\kappa_1}\le d\le \frac{Y}{\gamma}}
e\left(\frac{n\gamma}{\alpha_1}d\right)\right|.
\end{gather}
Обозначим $\tilde{\gamma}=\frac{\gamma}{(\gamma,\alpha_1)}$ и $\tilde{\alpha_1}=\frac{\alpha_1}{(\gamma,\alpha_1)}.$ Тогда
$\frac{n\gamma}{\alpha_1}=\frac{n\tilde{\gamma}}{\tilde{\alpha_1}},$ $(\tilde{\gamma},\tilde{\alpha_1})=1$ и
\begin{gather}\label{f-2-24}
\sum_{\frac{Y}{\kappa_1}\le d\le Y}\delta_{\alpha_1}(d-v')\delta_{\gamma}(d)\le
\frac{(\gamma,\alpha_1)}{\alpha_1}\sum_{n=1}^{\tilde{\alpha_1}}\left|\sum_{\frac{Y}{\gamma\kappa_1}\le d\le \frac{Y}{\gamma}}
e\left(\frac{n\tilde{\gamma}}{\tilde{\alpha_1}}d\right)\right|.
\end{gather}
Для дальнейших преобразований нам потребуется следующее утверждение, доказанное Н.М. Коробова ~\cite[гл. 1, \S1 лемма 1, лемма 3]{Korobov}. Пусть $Q$---целое и $P$---натуральное, пусть $q$---произвольное натуральное число, $1\le a<q$ и $(a,q)=1.$ Тогда справедлива оценка
\begin{gather}\label{f-2-24-1}
\sum_{n=1}^q\left|\sum_{x=Q+1}^{Q+P}e\left(\frac{na}{q}x\right)\right|
\le P+q\log q.
\end{gather}
Применяя \eqref{f-2-24-1} к \eqref{f-2-24}, получаем:
\begin{gather}\label{f-2-25}
\sum_{\frac{Y}{\kappa_1}\le d\le Y}\delta_{\alpha_1}(d-v')\delta_{\gamma}(d)\le
\frac{(\gamma,\alpha_1)}{\alpha_1}\left(\frac{Y}{\gamma}+\tilde{\alpha_1}\log\tilde{\alpha_1}\right).
\end{gather}
Докажем, что $Y\ge\gamma\tilde{\alpha_1}\log\tilde{\alpha_1}.$ Действительно, из формулы \eqref{f-2-21} заключаем, что $\gamma\le\frac{q}{\alpha_1}.$  Следовательно, используя очевидную оценку $\tilde{\alpha_1}\le q,$ получаем, что достаточно проверить справедливость неравенства: $Y\ge q \log q,$ которое выполнено по условию. Теперь формула \eqref{f-2-25} принимает следующий вид:
\begin{gather}\label{f-2-26}
\sum_{\frac{Y}{\kappa_1}\le d\le Y}\delta_{\alpha_1}(d-v')\delta_{\gamma}(d)\ll
\frac{Y(\gamma,\alpha_1)}{\gamma\alpha_1}.
\end{gather}
Подставляя \eqref{f-2-26} в \eqref{f-2-21}, получаем:
\begin{gather}\label{f-2-27}
\sum_{\frac{Y}{\kappa_1}\le d\le Y}\delta_{\alpha_1}(d-v')\left(\frac{q}{\alpha_1},d\right)^{1/2}\ll\frac{Y}{\alpha_1}
\sum_{\gamma|\frac{q}{\alpha_1}}\frac{(\gamma,\alpha_1)}{\gamma^{1/2}}.
\end{gather}
Применяя \eqref{f-2-27} в \eqref{f-2-20}, получаем:
\begin{gather}\label{f-2-28}
\Sigma_1\ll_{\epsilon}\frac{Y^{2+\epsilon}}{K^2}q^{-1/2+\epsilon}\alpha_1^{-1/2}
\sum_{\gamma|\frac{q}{\alpha_1}}\frac{(\gamma,\alpha_1)}{\gamma^{1/2}}.
\end{gather}
Оценим сумму из правой части \eqref{f-2-28}
\begin{gather}
\sum_{\gamma|\frac{q}{\alpha_1}}\frac{(\gamma,\alpha_1)}{\gamma^{1/2}}\le
\sum_{\beta|\alpha_1}\beta\sum_{\gamma|\frac{q}{\alpha_1}}\frac{1}{\gamma^{1/2}}\delta_{\beta}(\gamma)=
\sum_{\beta|(\frac{q}{\alpha_1},\alpha_1)}\beta^{1/2}\sigma_{-1/2}\left(\frac{q}{\beta\alpha_1}\right)\ll_{\epsilon} \alpha_1^{1/2}q^{\epsilon}.\label{f-2-29}
\end{gather}
Подставляя \eqref{f-2-29} в \eqref{f-2-28}, получаем:
\begin{gather}\label{f-2-30-1}
\Sigma_1\le\frac{Y^{2+\epsilon}}{K^2}q^{-1/2+\epsilon}.
\end{gather}
Лемма доказана.
\end{proof}
\end{Le}
  \item Оценим $\Sigma_2.$
Из определения $\Sigma_1$ и $\Sigma_2,$ а так же из формулы  \eqref{f-2-30} получаем:
\begin{gather}\label{f-2-31}
\Sigma_2\ll_{\epsilon}\frac{1}{T}\frac{Y^{2+\epsilon}}{K^2}q^{-1/2+\epsilon}.
\end{gather}
  \item Оценим $\Sigma_3.$ Подставляя в определение $\Sigma_3$ выражение для $R_{0}(kx)$ из \eqref{k-14}, получим:
\begin{gather}\label{f-2-32}
\Sigma_3=O\left(\frac{T}{q}\sum_{\frac{Y}{\kappa_1}\le d\le Y}
\delta_{\alpha_1}(d-v')\frac{1}{qs}\sum_{k=1}^q\left|\sum_{l=1}^q K_s(\pm q,dl,dkx)\right|\right)=
O\left(\Sigma_1\frac{TK^2}{Y}\right).
\end{gather}
Аналогично тому как была получена оценка суммы $\Sigma_1,$ получаем:
\begin{gather}\label{f-2-33}
\Sigma_3\ll_{\epsilon}TY^{1+\epsilon}q^{-1/2+\epsilon}.
\end{gather}
  \item Оценим $\Sigma_4.$ Подставляя в определение $\Sigma_4$ выражение для $R_{1}(kx)$ из \eqref{k-15}, получим:
\begin{gather}\label{f-2-34}
\Sigma_4=O\left(\frac{T}{q^3}\sum_{\frac{Y}{\kappa_1}\le d\le Y}\delta_{\alpha_1}(d-v')
\mathop{{\sum}'}_{-s/2<m\le s/2}\frac{1}{|m|}\sum_{k=1}^q\sum_{l=1}^q \left|K_s(\pm q,m+dl,dk)\right|\right).
\end{gather}
Применяя к $\left|K_s(\pm q,m+dl,dkx)\right|$ лемму \ref{lemma-k-1} и  используя замечание \ref{zam-k-1}, получим, что для $\forall\epsilon>0$ выполнено
\begin{gather}\label{f-2-35}
\left|K_s(\pm q,m+dl,dkx)\right|\ll_{\epsilon}
s^{\epsilon}(q,s)^{\epsilon}(q(m+dl),qdkx,(m+dl)dkx,s)^{1/2}\sqrt{s}.
\end{gather}
Ввиду того, что $s=qd,$ имеем
\begin{gather*}
(q(m+dl),qdkx,(m+dl)dkx,s)\le(qm+sl,skx,s)=(qm,qd)=q(m,d).
\end{gather*}
Тогда из \eqref{f-2-35} следует, что
\begin{gather}\label{f-2-36}
\left|K_s(\pm q,dl,dkx)\right|\ll_{\epsilon}
d^{\epsilon}q^{\epsilon}(m,d)^{1/2}\sqrt{sq}.
\end{gather}
Подставляя  \eqref{f-2-36} в \eqref{f-2-34}, получаем:
\begin{gather}\label{f-2-37}
\Sigma_4=O\left(\frac{T}{q^3}\sum_{\frac{Y}{\kappa_1}\le d\le Y}\delta_{\alpha_1}(d-v')q^2d^{1/2+\epsilon}q^{1+\epsilon}
\mathop{{\sum}'}_{-s/2<m\le s/2}\frac{(m,d)^{1/2}}{|m|}\right).
\end{gather}
Оценим сумму из правой части  \eqref{f-2-37}.
\begin{gather}\notag
\mathop{{\sum}'}_{-s/2<m\le s/2}\frac{(m,d)^{1/2}}{|m|}\ll\sum_{0<m\le s}\frac{(m,d)^{1/2}}{m}\le
\sum_{\beta|d}\beta^{1/2}\sum_{0<m\le s}\frac{1}{m}\delta_{\beta}(m)\le\\\le
\sum_{\beta|d}\beta^{-1/2}\log\frac{s}{\beta}\ll_{\epsilon} d^{\epsilon}q^{\epsilon},\label{f-2-38}
\end{gather}
для $\forall\epsilon>0.$ Применяя оценку  \eqref{f-2-38} к \eqref{f-2-37}, получим:
\begin{gather}\label{f-2-39}
\Sigma_4\ll_{\epsilon} Tq^{\epsilon}\sum_{\frac{Y}{\kappa_1}\le d\le Y}\delta_{\alpha_1}(d-v')d^{1/2+\epsilon}.
\end{gather}
Использование тривиальной оценки правой части дает
\begin{gather}\label{f-2-40}
\Sigma_4\ll_{\epsilon} Tq^{\epsilon}Y^{3/2+\epsilon}.
\end{gather}
  \item Оценим $\Sigma_5.$ Подставляя в определение $\Sigma_5$ выражение для $R_{2}(kx)$ из \eqref{k-16}, получим:
\begin{gather}\label{f-2-41}
\Sigma_5=O\left(\frac{T}{q^2}\sum_{\frac{Y}{\kappa_1}\le d\le Y}\delta_{\alpha_1}(d-v')
\mathop{{\sum}'}_{-s/2<n\le s/2}\frac{1}{|n|}\sum_{k=1}^q\sum_{l=1}^q \left|K_s(\pm q,dl,n+dkx)\right|\right).
\end{gather}
Совершенно аналогично тому, как была получена оценка \eqref{f-2-40}, поучаем:
\begin{gather}\label{f-2-42}
\Sigma_5\ll_{\epsilon} Tq^{1+\epsilon}Y^{3/2+\epsilon}.
\end{gather}
  \item Оценим $\Sigma_6.$ Подставляя в определение $\Sigma_6$ выражение для $R_{3}(kx)$ из \eqref{k-17}, получим:
\begin{gather}\label{f-2-43}
\Sigma_6=O\left(\frac{T}{q^2}\sum_{\frac{Y}{\kappa_1}\le d\le Y}\delta_{\alpha_1}(d-v')
\mathop{{\sum}'}_{-s/2<m\le s/2}\mathop{{\sum}'}_{-s/2<n\le s/2}\frac{1}{|mn|}
\sum_{k=1}^q\sum_{l=1}^q\left|K_s(\pm q,m+dl,n+dkx)\right|\right).
\end{gather}
Применяя к $\left|K_s(\pm q,m+dl,n+dkx)\right|$ лемму \ref{lemma-k-1} и  используя замечание \ref{zam-k-1}, получим, что для $\forall\epsilon>0$ выполнено
\begin{gather}\label{f-2-44}
\left|K_s(\pm q,m+dl,n+dkx)\right|\ll_{\epsilon}
d^{\epsilon}q^{\epsilon}(m,n,d)^{1/2}\sqrt{sq}.
\end{gather}
Подставляя  \eqref{f-2-44} в \eqref{f-2-43}, получаем:
\begin{gather}\label{f-2-45}
\Sigma_6\ll_{\epsilon} \frac{T}{q^2}\sum_{\frac{Y}{\kappa_1}\le d\le Y}\delta_{\alpha_1}(d-v')
q^{2}d^{1/2+\epsilon}q^{1+\epsilon}
\mathop{{\sum}'}_{-s/2<m\le s/2}\mathop{{\sum}'}_{-s/2<n\le s/2}\frac{(m,n,d)^{1/2}}{|mn|}.
\end{gather}
Оценим сумму из правой части  \eqref{f-2-45}.
\begin{gather}\notag
\mathop{{\sum}'}_{-s/2<m\le s/2}\mathop{{\sum}'}_{-s/2<n\le s/2}\frac{(m,n,d)^{1/2}}{|mn|}
\ll\sum_{0<m\le s}\sum_{0<n\le s}\frac{(m,n,d)^{1/2}}{mn}=\\\notag=
\sum_{\beta|d}\beta^{1/2}\sum_{0<m\le s}\sum_{0<n\le s}\frac{1}{mn}\delta_{\beta}((m,n))\le
\sum_{\beta|d}\beta^{1/2}\sum_{0<m\le s/\beta}\sum_{0<n\le s/\beta}\frac{1}{\beta^2mn}\le\\\le
\sum_{\beta|d}\beta^{-3/2}\log^2\frac{s}{\beta}\ll_{\epsilon}  s^{\epsilon},\label{f-2-46}
\end{gather}
для $\forall\epsilon>0.$ Подставляя  \eqref{f-2-46} в \eqref{f-2-45}, получаем
\begin{gather}\label{f-2-47}
\Sigma_6\ll_{\epsilon} Tq^{1+\epsilon}Y^{3/2+\epsilon}.
\end{gather}
\end{enumerate}
Тем самым мы полностью завершили доказательство теоремы.
\end{proof}
\end{Th}

\begin{Zam}\label{zam-k-3}
Используя суммы из замечания \ref{zam-k-2}, можно доказать, что при $Y>K^{4}q^6$ справедлива следующая оценка мощности множества $\M:$
\begin{gather}\label{f-3-43}
\left|\M\right|\le\frac{Y^{2+\epsilon}}{K^2}\frac{q^{\epsilon}}{q}(q,v)^{1/2}
\end{gather}
Однако, использование оценки \eqref{f-2-1} дает лучшую итоговую оценку на $\delta$ в теореме \ref{main}. Это объясняется зависимостью оценки \eqref{f-3-43} от вектора $\eta'=(u,v)^t.$
\end{Zam}
\section{Оценки тригонометрических сумм.}
Определим тригонометрическую сумму $S_N(\theta)$ следующим образом
\begin{gather}\label{def Sn 0}
S_N(\theta)=\sum_{\gamma\in\Omega_{N} }e(\theta\|\gamma\|),
\end{gather}
где $\Omega_{N}$~--особое множество матриц (ансамбль), построенное в  ~\cite[глава II]{FK}. В ~\cite[§7]{FK} было доказано, что для доказательства теоремы \ref{main} достаточно получить следующую оценку
\begin{gather}\label{result}
\int_0^1\left|S_N(\theta)\right|^2d\theta\ll\frac{1}{N}|\Omega_{N}|^2.
\end{gather}
В данном параграфе мы приведем результаты из ~\cite[§14]{FK}, необходимые для доказательства оценки \eqref{result}.\par
Из теоремы Дирихле следует, что для любого $\theta\in[0,1]$ найдутся $a,q\in\N\cup\{0\}$ и $\beta\in\rr$ такие что:
\begin{gather}\label{14-26}
\theta=\frac{a}{q}+\beta,\;(a,q)=1,\; 0\le a\le q\le N^{1/2},\;\beta=\frac{K}{N},\; |K|\le\frac{N^{1/2}}{q},
\end{gather}
причем $a=0$ или $a=q$ возможно только при $q=1.$ Обозначим
\begin{gather}\label{14-27}
P_{Q_1,Q}^{(\beta)}=\left\{
\theta=\frac{a}{q}+\beta\;\Bigl|\;(a,q)=1,\;  0\le a\le q,\;Q_1\le q\le Q
\right\}.
\end{gather}
Для каждого $q$ из $Q_1\le q\le Q$ определим каким-либо способом число $a_q,$ такое что $(a_q,q)=1,\,  0\le a_q\le q.$ Обозначим
\begin{gather}\label{15-0}
Z^{*}=\left\{
\theta=\frac{a_q}{q}+\beta\;\Bigl|\;Q_1\le q\le Q
\right\}.
\end{gather}

Положим также
\begin{gather}\label{15-8-0}
Q_0=\max\left\{\exp\left(\frac{10^5A^4}{\epsilon_0^2}\right),\exp(\epsilon_0^{-5})\right\},\quad \KK=\max\{1,|K|\}.
\end{gather}
\begin{Le}\label{lemma-15-2}
Если $\KK^2Q^{3}\le\frac{N^{1-\epsilon_0}}{12000A^4},\,\KK Q\ge Q_0,$ то имеет место оценка
\begin{gather}\label{15-9}
\sum_{\theta\in P_{Q_1,Q}^{(\beta)}}\left|S_N(\theta)\right|^2\ll
|\Omega_N|^2\KK^{12\epsilon_0}Q^{20\epsilon_0}\frac{\KK^{4(1-\delta)}Q^{6(1-\delta)+1}}{\KK Q^2_1}
\end{gather}
\end{Le}
\begin{Le}\label{lemma-15-2-1}
Если $\KK^2q^{2}\le\frac{N^{1-\epsilon_0}}{12000A^4},\,\KK q\ge Q_0,$ то имеет место оценка
\begin{gather}\label{15-18}
\left|S_N(\theta)\right|\ll
|\Omega_N|(\KK q)^{6\epsilon_0}\frac{(\KK q)^{2(1-\delta)}}{\KK^{1/2}q}.
\end{gather}
\end{Le}
\begin{Le}\label{lemma-15-3}
Если $\KK q\ge Q_0,$ то имеет место оценка
\begin{gather}\label{15-22}
|S_N(\theta)|\ll|\Omega_N|\frac{(\KK q)^{2\epsilon_0}N^{1-\delta+\epsilon_0}}{\KK q}.
\end{gather}
\end{Le}
\begin{Le}\label{lemma-15-4}
Пусть выполнены неравенства
$N^{\epsilon_0/2}\le Q^{1/2}\le Q_1\le Q,\,\KK Q\le N^{\alpha},$
причем $\frac{1}{4}<\alpha\le\frac{1}{2}+\epsilon_0,$ тогда имеет место оценка
\begin{gather}\label{15-27}
\sum_{\theta\in P_{Q_1,Q}^{(\beta)}}|S_N(\theta)|\ll
|\Omega_N|\left(
N^{1/2+\alpha/2-\delta+3\epsilon_0}Q+
N^{1-\delta+3\epsilon_0}\frac{Q^{1/2}}{\KK}\right).
\end{gather}
\end{Le}
Сформулируем еще одну лемму общего характера, доказанную в ~\cite[§12]{FK} Похожее утверждение использовалось С.В. Конягиным в ~\cite[следствие 17]{Konyagin}.
\begin{Le}\label{lemma-14-5}
Пусть $W$~--произвольное, конечное подмножество отрезка $[0,1]$ и $|W|>10.$ Пусть $f:W\rightarrow \rr_{+}$~--произвольная функция такая, что для любого подмножества $Z\subseteq W$ выполнена оценка
\begin{gather*}
\sum_{\theta\in Z}f(\theta)\le c_1|Z|^{1/2}+c_2,
\end{gather*}
где $c_1,c_2$~--неотрицательные константы, не зависящие от множества $Z.$ Тогда  справедлива оценка
\begin{gather}\label{14-49}
\sum_{\theta\in W}f^2(\theta)\ll c_1^2\log|W|+c_2\max_{\theta\in W}f(\theta)
\end{gather}
с абсолютной константой в знаке Виноградова.
\end{Le}

\section{«Случай $\mu=2$.»}
Данный параграф соответствует параграфу 15 в ~\cite{FK} и поэтому носит такое же название. Наша цель~--передоказать леммы из ~\cite[§15]{FK}, заменяя в них лемму \ref{BK lemma 7.1} на теорему \ref{theoremF-2} (нумерация утверждений соответствует настоящей работе).  Пусть $Z\subseteq P_{Q_1,Q}^{(\beta)}$ и $\KK=\max\{1,|K|\}.$ Для полноты изложения приведем результаты из ~\cite{FK}, которые понадобятся нам при доказательстве леммы \ref{lemma-16-2}.
\begin{Th}\label{theorem13-3-1}(~\cite[теорема 11.4.]{FK})
Пусть $M^{(2)}\ge(M^{(1)})^{2\epsilon_0},$ а неравенство \eqref{13-42}
\begin{gather}\label{13-42}
\exp\left(\frac{10^5A^4}{\epsilon_0^2}\right)\le M\le N\exp\left(-\frac{10^5A^4}{\epsilon_0^2}\right),
\end{gather}
выполнено при $M=M^{(1)}$ и $M=M^{(1)}M^{(2)}.$ Тогда ансамбль $\Omega_N$ может быть представлен в виде $\Omega_N=\Omega^{(1)}\Omega^{(2)}\Omega^{(3)},$ и для любых $\gamma_1\in\Omega^{(1)},\,\gamma_2\in\Omega^{(2)},\,\gamma_3\in\Omega^{(3)}$ справедливы неравенства
\begin{gather}\label{13-48-5}
\frac{M^{(2)}}{150A^2(M^{(1)})^{2\epsilon_0}}\le\|\gamma_2\|\le73A^2\frac{(M^{(2)})^{1+2\epsilon_0}}{(M^{(1)})^{2\epsilon_0}},
\end{gather}
\begin{gather}\label{13-48-6}
\frac{N}{150A^2(M^{(1)}M^{(2)})^{1+2\epsilon_0}}\le\|\gamma_3\|\le\frac{73A^2N}{M^{(1)}M^{(2)}}.
\end{gather}
\end{Th}
Будем записывать числа $\theta^{(1)},\theta^{(2)}\in P_{Q_1,Q}^{(\beta)}$ следующим образом:
\begin{gather}\label{14-28}
\theta^{(1)}=\frac{a^{(1)}}{q^{(1)}}+\beta,\quad \theta^{(2)}=\frac{a^{(2)}}{q^{(2)}}+\beta.
\end{gather}
Обозначим
\begin{gather}\label{14-42}
\M(g_{3})=\left\{ (g^{(1)}_{2},g^{(2)}_{2},\theta^{(1)},\theta^{(2)})\in
\Omega^{(2)}\times\Omega^{(2)}\times Z^2\Bigl|\,
\eqref{14-43}\, \mbox{и}\, \eqref{14-44}\, \mbox{выполнены}
\right\}
\end{gather}
где
\begin{gather}\label{14-43}
|g_2^{(1)}g_3-g_2^{(2)}g_3|_{1,2}\le \frac{73A^2N}{M^{(1)}\KK},
\end{gather}
\begin{gather}\label{14-44}
\|g_2^{(1)}g_3\frac{a^{(1)}}{q^{(1)}}-g_2^{(2)}g_3\frac{a^{(2)}}{q^{(2)}}\|_{1,2}=0.
\end{gather}
\begin{Le}\label{lemma-14-4}
Пусть выполнены условия  теоремы \ref{theorem13-3-1} и на ее основе построено разложение $\Omega_N$ в виде $\Omega_N=\Omega^{(1)}\Omega^{(2)}\Omega^{(3)},$. Пусть $M^{(1)}$ такое, что для любых $\theta^{(1)},\theta^{(2)}\in Z$ выполнено
\begin{gather}\label{14-45}
[q^{(1)},q^{(2)}]<\frac{M^{(1)}}{74A^2\KK}.
\end{gather}
Тогда имеет место оценка
\begin{gather}\label{14-46}
\sum_{\theta\in Z}\left|S_N(\theta)\right|\ll
(M^{(1)})^{1+2\epsilon_0}\left|\Omega^{(1)}\right|^{1/2}\sum_{g_{3}\in\Om^{(3)}}
\left|\M(g_{3})\right|^{1/2}.
\end{gather}
где $\Om^{(3)}=\Omega^{(3)}(0,1)^{t}.$
 \end{Le}

\begin{Le}\label{lemma-16-2}
Если, выполнены неравенства
\begin{gather}\label{16-13-00}
\KK^{5+28\epsilon_0} Q^{5+21\epsilon_0}<N^{1-\epsilon_0},\quad \KK Q\ge Q_0,
\end{gather}
то имеет место оценка
\begin{gather}\label{16-13}
\sum_{\theta\in P_{Q_1,Q}^{(\beta)}}\left|S_N(\theta)\right|^2\ll|\Omega_N|^2
\frac{\KK^{10(1-\delta)} Q^{10(1-\delta)+1}\KK^{60\epsilon_0}Q^{60\epsilon_0}}{\KK^2Q_1^{1/2}}.
\end{gather}
\begin{proof}
Доказательство леммы \ref{lemma-16-2} практически не отличается от доказательства леммы 15.1. из ~\cite{FK}, за исключением замены леммы \ref{BK lemma 7.1} на теорему \ref{theoremF-2} и другого выбора параметров $M^{(1)},M^{(2)}.$ Используя тривиальную оценку
\begin{gather}\label{15-16-1}
\sum_{\theta\in P_{Q_1,Q}^{(\beta)}}\left|S_N(\theta)\right|^2\le
Q\sum_{Q_1\le q\le Q}\max\limits_{1\le a\le q, (a,q)=1}\left|S_N(\frac{a}{q}+\frac{K}{N})\right|^2=Q\sum_{\theta\in Z^*}\left|S_N(\theta)\right|^2,
\end{gather}
где в качестве $a_q$ выбраны числители, на которых достигается максимум, получаем, что достаточно доказать неравенство
\begin{gather}\label{16-1-2}
\sum_{\theta\in Z^{*}}\left|S_N(\theta)\right|^2\ll|\Omega_N|^2
\frac{\KK^{10(1-\delta)} Q^{10(1-\delta)}\KK^{60\epsilon_0}Q^{60\epsilon_0}}{\KK^2Q_1^{1/2}}.
\end{gather}
Пусть $Z\subseteq Z^{*}$~--любое подмножество. Воспользуемся разложением ансамбля $\Omega_N$ по  по теореме \ref{theorem13-3-1}. Положим:
\begin{gather}\label{16-14}
M^{(1)}=75A^2\KK Q^2,\quad M^{(2)}=(75A^2\KK^4Q^3)^{1+7\epsilon_0},
\end{gather}
тогда выполнено условие \eqref{14-45} и все условия теоремы \ref{theorem13-3-1}. Следовательно, выполнены все условия леммы \ref{lemma-14-4} и значит имеет место оценка \eqref{14-46}. Аналогично тому, как это было сделано в ~\cite[лемме 14.1.]{FK}, из \eqref{14-44} получаем: $q^{(1)}=q^{(2)}=\q$ и, следовательно, $a^{(1)}=a^{(2)}.$ Тогда соотношения  \eqref{14-43} и \eqref{14-44} дают:
\begin{gather}\label{16-4}
\left|(g_2^{(1)}-g_2^{(2)})g_3\right|_{1,2}\le \frac{73A^2N}{M^{(1)}\KK},\quad
\left((g_2^{(1)}-g_2^{(2)})g_3\right)_{1,2}\equiv 0\pmod{\q},
\end{gather}
Положим:
\begin{gather}\label{16-5}
\eta'=g_2^{(2)}g_3,\,\eta=g_3,\,\gamma=g_2^{(1)},\,X=\|\eta'\|,\,Y=\|g_2^{(2)}\|,\,K_1=\KK\frac{XM^{(1)}}{73A^2N}.
\end{gather}
Без ограничения общности можно считать, что $\|g_2^{(1)}\|\le\|g_2^{(2)}\|,$ тогда  $\|\gamma\|\le Y.$ Из свойств ансамбля так же следует, что $\|\gamma\|\asymp Y.$ Кроме того, из теоремы \ref{theorem13-3-1} следует, что
\begin{gather}\label{16-5-1}
\frac{\KK}{(M^{(1)})^{4\epsilon_0}(M^{(2)})^{2\epsilon_0}}\ll K_1\ll\KK\frac{(M^{(2)})^{2\epsilon_0}}{(M^{(1)})^{2\epsilon_0}},\quad
\frac{(M^{(2)})}{(M^{(1)})^{2\epsilon_0}}\ll Y\ll\frac{(M^{(2)})^{1+2\epsilon_0}}{(M^{(1)})^{2\epsilon_0}},\,
\end{gather}
Соотношения  \eqref{16-4} могут быть записаны в следующем виде:
\begin{gather}\label{16-7}
|\gamma\eta-\eta'|_{1,2}<\frac{X}{K_1},\quad(\gamma\eta-\eta')_{1,2}\equiv0\pmod{q}.
\end{gather}
Проверим выполнены ли условия теоремы \ref{theoremF-2}. Для этого достаточно убедиться, что $Y<X,\,K_1^{4}Q^3<Y,$ то есть ~--проверить, что
\begin{gather}\label{16-8-1}
\KK^4 Q^{3}<(M^{(2)})^{1-8\epsilon_0}(M^{(1)})^{6\epsilon_0}.
\end{gather}
Неравенства \eqref{16-8-1} выполнены в силу выбора параметров $M^{(1)},M^{(2)}.$ Таким образом для оценки мощности множества $\M(g_3)$ может быть применена теорема \ref{theoremF-2} следующим образом. Фиксируем $g_2^{(2)}$ любым из $|\Omega^{(2)}|$ способов. Так же фиксируем $\frac{a^{(1)}}{q^{(1)}}$ любым из $|Z|$ способов. По доказанному, этим определено $\frac{a^{(2)}}{q^{(2)}}.$ Следовательно, при фиксированном $g_3$ получаем:
\begin{gather*}\label{16-17}
|\M(g_3)|\ll\frac{Y^2(YQ)^{\epsilon_0}}{K_1^2Q_1^{1/2}}\left|\Omega^{(2)}\right||Z|
\end{gather*}
Используя оценки \eqref{16-5-1}, получаем:
\begin{gather}\label{16-18}
|\M(g_3)|\ll\frac{(M^{(2)})^{2+8\epsilon_0}(M^{(1)})^{4\epsilon_0}(\KK Q)^{5\epsilon_0}}{\KK^2 Q_1^{1/2}}
\left|\Omega^{(2)}\right||Z|\ll
\frac{\KK^8Q^6\KK^{98\epsilon_0}Q^{80\epsilon_0}}{\KK^2Q_1^{1/2}}\left|\Omega^{(2)}\right||Z|.
\end{gather}
Подставляя \eqref{16-18} в  \eqref{14-46}, получаем:
\begin{gather}\label{16-20}
\sum_{\theta\in Z}\left|S_N(\theta)\right|\ll|Z|^{1/2}
(M^{(1)})^{1+2\epsilon_0}\left|\Omega^{(1)}\right|^{1/2}\left|\Omega^{(3)}\right|\left|\Omega^{(2)}\right|^{1/2}
\frac{\KK^4Q^3\KK^{49\epsilon_0}Q^{40\epsilon_0}}{\KK Q_1^{1/4}}.
\end{gather}
Используя оценки  $\left|\Omega^{(1)}\right|\ge(M^{(1)})^{\delta-2\epsilon_0},\,\left|\Omega^{(2)}\right|\ge(M^{(2)})^{\delta-2\epsilon_0},$
доказанные в ~\cite[§11]{FK}, получаем:
\begin{gather}\label{16-21}
\sum_{\theta\in Z}\left|S_N(\theta)\right|\ll|Z|^{1/2}
|\Omega_N|\frac{(M^{(1)})^{1-\delta+2,5\epsilon_0}}{(M^{(2)})^{\delta-\epsilon_0/2}}
\frac{\KK^4Q^3\KK^{49\epsilon_0}Q^{40\epsilon_0}}{\KK Q_1^{1/4}}.
\end{gather}
Подставляя \eqref{16-14} в  \eqref{16-21}, получаем:
\begin{gather}\label{16-22}
\sum_{\theta\in Z}\left|S_N(\theta)\right|\ll|Z|^{1/2}|\Omega_N|
\frac{\KK^{5(1-\delta)}Q^{5(1-\delta)}\KK^{30\epsilon_0}Q^{30\epsilon_0}}{\KK Q_1^{1/4}}.
\end{gather}
Применяя лемму \ref{lemma-14-5} с $W=Z^*,c_2=0,f(\theta)=\frac{|S_N(\theta)|}{|\Omega_N|},$ получаем \eqref{16-1-2}. Лемма доказана.
\end{proof}
\end{Le}
\begin{Le}\label{lemma-16-2-0}
Если выполнено неравенство  $\KK q>Q_0,$ то имеет место оценка
\begin{gather}\label{16-13-0}
|S_N(\theta)|\ll_{A,\epsilon_0}|\Omega_N|
\frac{\KK^{5(1-\delta)}q^{4(1-\delta)}\KK^{30\epsilon_0}q^{24\epsilon_0}}{\KK q^{1/4}}.
\end{gather}
\begin{proof}
Положим  $Z=\{\theta\}$ и
\begin{gather}\label{16-30}
M^{(1)}=75A^2\KK q,\quad M^{(2)}=(75A^2\KK^4q^3)^{1+7\epsilon_0}.
\end{gather}
тогда выполнено условие \eqref{14-45} из леммы \ref{lemma-14-4}. Предположим, что справедливо неравенство
\begin{gather}\label{16-31}
\KK^{5+34\epsilon_0}q^{4+26\epsilon_0}<N,
\end{gather}
тогда выполнены все условия теоремы \ref{theorem13-3-1}, а, следовательно, и все условия леммы \ref{lemma-14-4}. Значит имеет место оценка \eqref{14-46}. Далее мы будем использовать обозначения \eqref{16-5}. Проверка выполнимости условий теоремы \ref{theoremF-2} проводится аналогично. Тем самым для оценки мощности множества $\M(g_3)$ может быть применена теорема \ref{theoremF-2}:
\begin{gather}\label{16-32}
|\M(g_3)|\ll\frac{(M^{(2)})^{2+8\epsilon_0}(M^{(1)})^{4\epsilon_0}(\KK q)^{5\epsilon_0}}{\KK^2 q^{1/2}}
\left|\Omega^{(2)}\right||\ll
\KK^6q^{5,5}\KK^{98\epsilon_0}q^{76\epsilon_0}\left|\Omega^{(2)}\right|.
\end{gather}
Подставляя \eqref{16-32} в  \eqref{14-46}, получаем
\begin{gather}\label{16-33}
\left|S_N(\theta)\right|\ll
(M^{(1)})^{1+2\epsilon_0}\left|\Omega^{(1)}\right|^{1/2}\left|\Omega^{(3)}\right|\left|\Omega^{(2)}\right|^{1/2}
\KK^3q^{11/4}\KK^{49\epsilon_0}q^{38\epsilon_0}
\end{gather}
Аналогично тому как были получены формулы \eqref{16-21} и \eqref{16-22} приходим к оценке
\begin{gather}\label{16-34}
\left|S_N(\theta)\right|\ll_{A}\Omega_N|
\frac{\KK^{5(1-\delta)}q^{4(1-\delta)}\KK^{30\epsilon_0}q^{24\epsilon_0}}{\KK q^{1/4}}.
\end{gather}
Пусть теперь неравенство \eqref{16-31} не выполнено, то есть
\begin{gather}\label{16-35}
N\le\KK^{5+34\epsilon_0}q^{4+26\epsilon_0},
\end{gather}
тогда выполнены условия леммы \ref{lemma-15-3} и учитывая \eqref{16-35}, получаем
\begin{gather}\label{16-36}
|S_N(\theta)|\ll
|\Omega_N|\frac{\KK^{5(1-\delta)}q^{4(1-\delta)}\KK^{30\epsilon_0}q^{24\epsilon_0}}{\KK q}.
\end{gather}
Лемма доказана.
\end{proof}
\end{Le}

\section{Оценка интегралов от $|S_N(\theta)|^2$.}
Доказательство теоремы \ref{main} схоже с доказательством теоремы \ref{uslov} в ~\cite{FK}. Нам понадобится ряд лемм из ~\cite[§16]{FK}, которые мы приведем без доказательства.
\begin{Le}\label{lemma-17-1}
Имеет место неравенство
\begin{gather}\label{17-1}
\int_0^1\left|S_N(\theta)\right|^2d\theta\le \frac{1}{N}
\mathop{{\sum}^*}_{0\le a\le q\le N^{1/2}}\int\limits_{|K|\le\frac{N^{1/2}}{q}}
\left|S_N(\frac{a}{q}+\frac{K}{N})\right|^2dK,
\end{gather}
где $\mathop{{\sum}^*}$ означает сумму по взаимно простым $a$ и $q$ при $q\ge1,$ и $a=0,1$ при $q=1.$
\end{Le}
Напомним, что
\begin{gather*}
Q_0=\max\left\{\exp\left(\frac{10^5A^4}{\epsilon_0^2}\right),\exp(\epsilon_0^{-5})\right\}.
\end{gather*}
\begin{Le}\label{lemma-17-2}
Имеет место неравенство
\begin{gather}\notag
\int_0^1\left|S_N(\theta)\right|^2d\theta\le
2Q_0^2\frac{|\Omega_N|^2}{N}+
\frac{1}{N}\mathop{{\sum}^*}_{0\le a\le q\le N^{1/2}\atop q>Q_0}\int\limits_{\frac{Q_0}{q}\le|K|\le\frac{N^{1/2}}{q}}
\left|S_N(\frac{a}{q}+\frac{K}{N})\right|^2dK+\\
\frac{1}{N}\mathop{{\sum}^*}_{0\le a\le q\le Q_0}\int\limits_{\frac{Q_0}{q}\le|K|\le\frac{N^{1/2}}{q}}
\left|S_N(\frac{a}{q}+\frac{K}{N})\right|^2dK
+\frac{1}{N}\mathop{{\sum}^*}_{1\le a\le q\le N^{1/2}\atop q>Q_0}\int\limits_{|K|\le\frac{Q_0}{q}}
\left|S_N(\frac{a}{q}+\frac{K}{N})\right|^2dK\label{17-4}
\end{gather}
\end{Le}
Оценим сначала третий из интегралов в правой части \eqref{17-4}.  Нам будет удобно далее использовать следующие обозначения:
\begin{gather}\label{17-7}
\gamma=1-\delta,\quad\xi_1=N^{2\gamma+6\epsilon_0}.
\end{gather}
Из ~\cite[лемма 16.2, лемма 16.3]{FK} следует следующее утверждение.
\begin{Le}\label{lemma-17-3-0}
При $\gamma<\frac{1}{8}$ и $\epsilon_0\in(0,\frac{1}{1000})$  имеет место неравенство
\begin{gather}\label{17-7-1}
\frac{1}{N}\mathop{{\sum}^*}_{1\le a\le q\le N^{1/2}}\int\limits_{|K|\le\frac{Q_0}{q}}
\left|S_N(\frac{a}{q}+\frac{K}{N})\right|^2dK\ll\frac{|\Omega_N|^2}{N}.
\end{gather}
\end{Le}
Оценим второй из интегралов в правой части \eqref{17-4}.
\begin{Le}\label{lemma-17-3-2}
При $\gamma<\frac{1}{12}$ и $\epsilon_0\in(0,\frac{1}{1000})$ имеет место неравенство
\begin{gather}\label{17-7-8}
\frac{1}{N}\mathop{{\sum}^*}_{0\le a\le q\le Q_0}\int\limits_{\frac{Q_0}{q}\le|K|\le\frac{N^{1/2}}{q}}
\left|S_N(\frac{a}{q}+\frac{K}{N})\right|^2dK\ll\frac{|\Omega_N|^2}{N}.
\end{gather}
\begin{proof}
Обозначим через $I$ интеграл в левой части \eqref{17-7-8}.  Применяя лемму \ref{lemma-16-2-0}, получаем
\begin{gather}\label{17-7-9}
I\ll|\Omega_N|^2q^{8\gamma-\frac{1}{2}+48\epsilon_0}
\int\limits_{\frac{Q_0}{q}\le|K|\le\frac{N^{1/2}}{q}}
\KK^{10\gamma-2+60\epsilon_0}dK\ll|\Omega_N|^2q^{\frac{1}{6}+48\epsilon_0},
\end{gather}
в силу выбора параметра $\gamma.$ Суммируя \eqref{17-7-9} по $0\le a\le q\le Q_0,$ получаем \eqref{17-7-8}. Лемма доказана.
\end{proof}
\end{Le}
Нам осталось оценить первый интеграл из правой части \eqref{17-4}, то есть
\begin{gather}\label{18-7-10}
\frac{1}{N}\mathop{{\sum}^*}_{0\le a\le q\le N^{1/2}\atop q>Q_0}\int\limits_{\frac{Q_0}{q}\le|K|\le\frac{N^{1/2}}{q}}
\left|S_N(\frac{a}{q}+\frac{K}{N})\right|^2dK
\end{gather}
Этому будут посвящены последующие леммы. Обозначим
\begin{gather}\label{18-1}
\xi_3=N^{4\gamma+15\epsilon_0}.
\end{gather}
Мы разобьем область суммирования и интегрирования по $q,\,K$ на пять подобластей см рис.

\begin{center}
  \includegraphics[width=450pt,height=350pt]{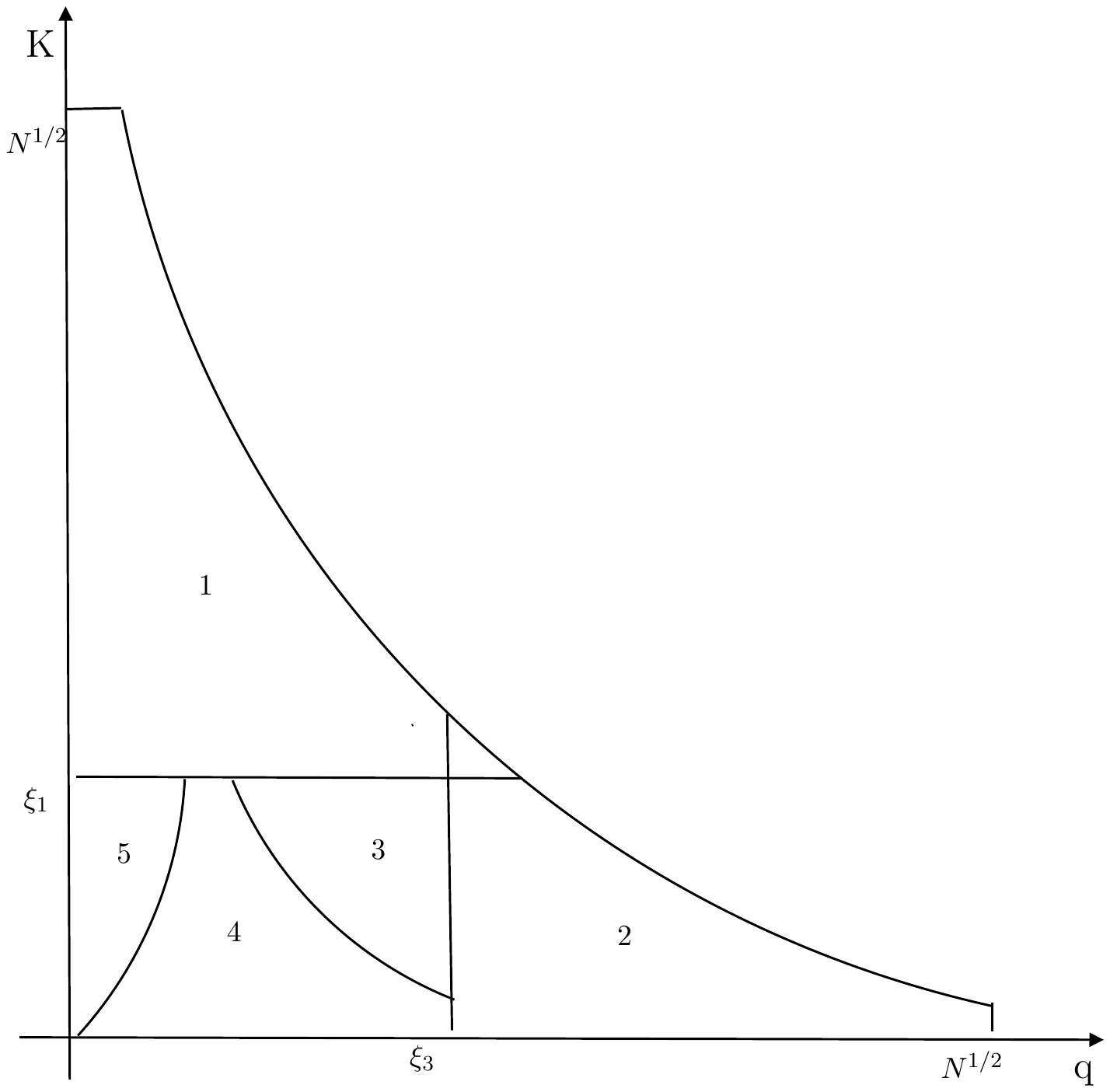}\\
\end{center}
области 1 соответствует лемма \ref{lemma-17-3}, области 2~--лемма \ref{lemma-17-7}, области 3~--лемма \ref{lemma-18-1}, области 4~--лемма \ref{lemma-18-2}, области 5~--лемма \ref{lemma-18-3}.

Леммы \ref{lemma-17-3}~--\ref{lemma-17-7}, следующие далее, также были доказаны в ~\cite[§16]{FK}.
\begin{Le}\label{lemma-17-3}
При $\gamma\le\frac{27-\sqrt{633}}{16}-5\epsilon_0,\,\epsilon_0\in(0,\frac{1}{2500})$ имеет место неравенство
\begin{gather}\label{17-7-11}
\frac{1}{N}\mathop{{\sum}^*}_{1\le a\le q\le N^{1/2}\atop q>Q_0}\int\limits_{\xi_1\le|K|\le\frac{N^{1/2}}{q}}
\left|S_N(\frac{a}{q}+\frac{K}{N})\right|^2dK\ll\frac{|\Omega_N|^2}{N}.
\end{gather}
\end{Le}
Пусть
\begin{gather*}
c_1=c_1(N),\,c_2=c_2(N),\,Q_0\le c_1<c_2\le  N^{1/2},
\end{gather*}
и пусть
\begin{gather*}
f_1=f_1(N,q),\,f_2=f_2(N,q),\,\frac{Q_0}{q}\le f_1<f_2\le  \frac{N^{1/2}}{q},\\
m_1=\min\{f_1(N,N_j),f_1(N,N_{j+1})\},\,m_2=\max\{f_2(N,N_j),f_2(N,N_{j+1})\}.
\end{gather*}
\begin{Le}\label{lemma-17-5}
Если функции $f_1(N,q),f_2(N,q)$ монотонны по $q,$ то имеет место неравенство
\begin{gather}\notag
\mathop{{\sum}^*}_{c_1\le q\le c_2\atop 1\le a\le q }\,\int\limits_{f_1\le|K|\le f_2}
\left|S_N(\frac{a}{q}+\frac{K}{N})\right|^2dK\le\\\le
\sum_{j:\,c_1^{1-\epsilon_0}\le N_j\le c_2\,}\int\limits_{m_1\le|K|\le m_2}
\mathop{{\sum}^*}_{N_j\le q\le N_{j+1}\atop 1\le a\le q }\left|S_N(\frac{a}{q}+\frac{K}{N})\right|^2dK.\label{17-7-15}
\end{gather}
\end{Le}
\begin{Le}\label{lemma-17-7}
При $\gamma\le\frac{27-\sqrt{633}}{16}-5\epsilon_0,\,\epsilon_0\in(0,\frac{1}{2500})$  имеет место неравенство
\begin{gather}\label{17-12}
\frac{1}{N}\mathop{{\sum}^*}_{1\le a\le q\le N^{1/2}\atop q>N^{4\gamma+14\epsilon_0}}
\int\limits_{\frac{Q_0}{q}\le|K|\le\frac{N^{1/2}}{q}}
\left|S_N(\frac{a}{q}+\frac{K}{N})\right|^2dK\ll\frac{|\Omega_N|^2}{N}.
\end{gather}
\end{Le}
Во всех последующих леммах для упрощения записи обозначим $Q_1=N_j,Q=N_{j+1}.$ Используя соотношения из ~\cite[лемма 9.1.]{FK}, получаем:
\begin{gather*}
\frac{Q}{Q_1}\le Q^{\epsilon_0}\le N^{\epsilon_0/2},\quad c_1^{1-\epsilon_0}\le Q_1\le c_2,\quad c_1\le Q\le c_2^{1+2\epsilon_0}.
\end{gather*}
Мы будем далее использовать эти оценки без ссылок на них. Напомним, что ${\KK=\max\{1,|K|\}.}$ Заметим, что у нас всегда будет выполнено $m_2\ge1,$ поэтому при $\eta<1$ имеем:
\begin{gather}\label{int-K-1}
\int\limits_{m_1\le|K|\le m_2}\frac{dK}{\KK^{\eta}}\ll m_2^{1-\eta}.
\end{gather}
При $\eta>1$ всегда выполнено:
\begin{gather}\label{int-K-2}
\int\limits_{m_1\le|K|\le m_2}\frac{dK}{\KK^{\eta}}<\frac{1}{m_1^{\eta-1}}.
\end{gather}
Однако, если $m_1\le1,$ то при $\eta>1$ имеем:
\begin{gather}\label{int-K-3}
\int\limits_{m_1\le|K|\le m_2}\frac{dK}{\KK^{\eta}}\ll1.
\end{gather}

\begin{Le}\label{lemma-18-1}
При $\gamma\le\frac{1}{12}-4\epsilon_0$ и $\epsilon_0\in(0,\frac{1}{1000})$  имеет место неравенство
\begin{gather}\label{18-2}
\frac{1}{N}\mathop{{\sum}^*}_{1\le a\le q\le \xi_3\atop q>\xi_3/\xi_1}
\int\limits_{\frac{\xi_3}{q}\le|K|\le\xi_1}
\left|S_N(\frac{a}{q}+\frac{K}{N})\right|^2dK\ll\frac{|\Omega_N|^2}{N}.
\end{gather}
\begin{proof}
Воспользуемся леммой \ref{lemma-17-5} с
\begin{gather*}
c_1=N^{2\gamma+9\epsilon_0},\,c_2=\xi_3,
f_1=\frac{\xi_3}{q},\,f_2=\xi_1,\,
m_1=\frac{\xi_3}{Q},\,m_2=\xi_1.
\end{gather*}
Оценивая $|S_N(\frac{a}{q}+\frac{K}{N})|$ через максимум по $a,q,$ получаем
\begin{gather}\label{17-9}
\mathop{{\sum}^*}_{N_j\le q\le N_{j+1}\atop 1\le a\le q }\left|S_N(\frac{a}{q}+\frac{K}{N})\right|^2\le
\max\limits_{Q_1\le q\le Q\atop 1\le a\le q, (a,q)=1}\left|S_N(\frac{a}{q}+\frac{K}{N})\right|
\mathop{{\sum}^*}_{Q_1\le q\le Q\atop 1\le a\le q }\left|S_N(\frac{a}{q}+\frac{K}{N})\right|.
\end{gather}
Заметим, что $KQ\le\xi_1^{1+2\epsilon_0}\xi_3\le N^{1/2},$ поэтому, применяя лемму \ref{lemma-15-3} для оценки максимума и лемму \ref{lemma-15-4} с $\alpha=\frac{1}{2},$ получаем:
\begin{gather}\label{17-10}
\mathop{{\sum}^*}_{N_j\le q\le N_{j+1}\atop 1\le a\le q }\left|S_N(\frac{a}{q}+\frac{K}{N})\right|^2\ll
|\Omega_N|^2\left(
\frac{N^{2\gamma-1/4+6\epsilon_0}Q}{(\KK Q_1)^{1-2\epsilon_0}}+
N^{2\gamma+5\epsilon_0}\frac{Q^{1/2}}{\KK^{2-2\epsilon_0}Q_1^{1-2\epsilon_0}}\right).
\end{gather}
Используя \eqref{int-K-1} при интегрировании по $K$ первого слагаемого и \eqref{int-K-2} при интегрировании второго, получаем:
\begin{gather*}
\int\limits_{m_1\le|K|\le m_2}\mathop{{\sum}^*}_{N_j\le q\le N_{j+1}\atop 1\le a\le q }
\left|S_N(\frac{a}{q}+\frac{K}{N})\right|^2dK\ll
|\Omega_N|^2\left(
N^{2\gamma-1/4+8\epsilon_0}\frac{Q}{Q_1}+
\frac{N^{2\gamma+6\epsilon_0}}{\xi_3^{1-2\epsilon_0}}\xi_3^{1/2}\right).
\end{gather*}
подставляя значения $\xi_3$ из \eqref{18-1} и используя  условие на $\gamma,$ получаем
\begin{gather}\label{17-11}
\int\limits_{m_1\le|K|\le m_2}\mathop{{\sum}^*}_{N_j\le q\le N_{j+1}\atop 1\le a\le q }
\left|S_N(\frac{a}{q}+\frac{K}{N})\right|^2dK\ll
|\Omega_N|^2 N^{-0,1\epsilon_0}.
\end{gather}
Поскольку количество слагаемых в сумме по $j$ в \eqref{17-7-15} превосходит $\log\log N,$ то
\begin{gather*}
\sum_{j:\,c_1^{1-\epsilon_0}\le N_j\le c_2\,}\int\limits_{m_1\le|K|\le m_2}
\mathop{{\sum}^*}_{N_j\le q\le N_{j+1}\atop 1\le a\le q }\left|S_N(\frac{a}{q}+\frac{K}{N})\right|^2dK\ll
|\Omega_N|^2 .
\end{gather*}
Лемма доказана.
\end{proof}
\end{Le}
Обозначим
\begin{gather}\label{18-3}
\nu=\frac{3+\sqrt{61}}{4}=2,70\ldots.
\end{gather}
\begin{Le}\label{lemma-18-2}
При $\gamma\le\frac{1}{9+\sqrt{61}}-4\epsilon_0$ и $\epsilon_0\in(0,\frac{1}{1000})$ имеет место неравенство
\begin{gather}\label{18-4}
\frac{1}{N}\mathop{{\sum}^*}_{1\le a\le q\le \xi_3\atop q>Q_0}
\int\limits_{\frac{Q_0}{q}\le|K|\le\min\{q^{\nu},\frac{\xi_3}{q}\}}
\left|S_N(\frac{a}{q}+\frac{K}{N})\right|^2dK\ll\frac{|\Omega_N|^2}{N}.
\end{gather}
\begin{proof}
Эта лемма доказывается совершенно аналогично лемме 16.12 из ~\cite[§16]{FK}. Необходимо лишь проверить условия леммы \ref{lemma-15-2} (в нумерации данной статьи), которая использовалась в доказательстве леммы 16.12 из ~\cite[§16]{FK}
\begin{gather*}
\KK^2Q^3\le\xi_3^3\le N^{1-1,5\epsilon_0}.
\end{gather*}
Лемма доказана.
\end{proof}
\end{Le}
\begin{Le}\label{lemma-18-3}
При $\gamma\le\frac{1}{9+\sqrt{61}}-8\epsilon_0$ и $\epsilon_0\in(0,\frac{1}{1000})$ имеет место неравенство
\begin{gather}\label{18-5}
\frac{1}{N}\mathop{{\sum}^*}_{1\le a\le q\le \xi_1^{1/\nu}\atop q>Q_0}
\int\limits_{q^{\nu}\le|K|\le\xi_1}
\left|S_N(\frac{a}{q}+\frac{K}{N})\right|^2dK\ll\frac{|\Omega_N|^2}{N}.
\end{gather}
\begin{proof}
Воспользуемся леммой \ref{lemma-17-5} с
\begin{gather*}
c_1=Q_0,\,c_2=\xi_1^{1/\nu},\,
f_1=q^{\nu},\,f_2=\xi_1,\,
m_1=Q_1^{\nu},\,m_2=\xi_1.
\end{gather*}
Проверим, что выполнены условия леммы \ref{lemma-16-2}. Для этого достаточно проверить, что $\KK^{5+35\epsilon_0}Q^{5+35\epsilon_0}<N,$ то есть, что
\begin{gather*}
(\xi_1\xi_1^{1/\nu})^{5+40\epsilon_0}<N,
\end{gather*}
Последнее неравенство выполнено ввиду условий на $\gamma,\epsilon_0.$ Применяя лемму \ref{lemma-16-2} и учитывая, что $\frac{Q}{Q_1}\le Q^{\epsilon_0},$ получаем:
\begin{gather}\label{18-6}
\mathop{{\sum}^*}_{N_j\le q\le N_{j+1}\atop 1\le a\le q }\left|S_N(\frac{a}{q}+\frac{K}{N})\right|^2\ll
|\Omega_N|^2
\KK^{10\gamma-2+60\epsilon_0}Q^{10\gamma+\frac{1}{2}+61\epsilon_0}.
\end{gather}
Используя \eqref{int-K-2} при интегрировании по $K,$ получаем:
\begin{gather*}
\int\limits_{m_1\le|K|\le m_2}
\mathop{{\sum}^*}_{N_j\le q\le N_{j+1}\atop 1\le a\le q }\left|S_N(\frac{a}{q}+\frac{K}{N})\right|^2dK\ll
|\Omega_N|^2Q^{10\gamma+\frac{1}{2}+61\epsilon_0}
Q_1^{\nu(10\gamma-1+60\epsilon_0)}\le\\\le
|\Omega_N|^2Q^{10\gamma+\frac{1}{2}+61\epsilon_0+\nu(10\gamma-1+60\epsilon_0)+3\epsilon_0}.
\end{gather*}
Для того чтобы сумма по $j$ была ограничена константой достаточно потребовать
\begin{gather*}
10(1+\nu)\gamma+\frac{1}{2}-\nu+230\epsilon_0\le-0,1\epsilon_0\,\Rightarrow\,
\gamma\le\frac{\nu-1/2}{10(1+\nu)}-8\epsilon_0=\frac{1}{9+\sqrt{61}}-8\epsilon_0
\end{gather*}
Последнее же неравенство выполнено ввиду условий на $\gamma.$ Лемма доказана.
\end{proof}
\end{Le}
\section{Доказательство теоремы \ref{main}.}
Пусть $\gamma<\frac{1}{9+\sqrt{61}}.$ Выбираем $\epsilon_0$ так, чтобы $\epsilon_0\in(0,\frac{1}{1000})$ и
$\gamma\le\frac{1}{9+\sqrt{61}}-8\epsilon_0.$ Тогда из лемм \ref{lemma-17-3}~--\ref{lemma-18-3} следует, что и первый интеграл из правой части \eqref{17-4} меньше чем $\frac{|\Omega_N|^2}{N},$ то есть
\begin{gather}\label{17-30}
\frac{1}{N}\mathop{{\sum}^*}_{0\le a\le q\le N^{1/2}\atop q>Q_0}\int\limits_{\frac{Q_0}{q}\le|K|\le\frac{N^{1/2}}{q}}
\left|S_N(\frac{a}{q}+\frac{K}{N})\right|^2dK\ll\frac{|\Omega_N|^2}{N}.
\end{gather}.
Подставляя \eqref{17-30} и результаты  лемм \ref{lemma-17-3-0}~--\ref{lemma-17-3-2} в лемму \ref{lemma-17-2}, получаем
\begin{gather}\label{17-31}
\int_0^1\left|S_N(\theta)\right|^2d\theta\ll\frac{|\Omega_N|^2}{N}\quad\mbox{при}\quad
\gamma\le\frac{1}{9+\sqrt{61}}-8\epsilon_0,\,\epsilon_0\in(0,\frac{1}{1000})
\end{gather}
Таким образом нами доказано неравенство \eqref{result}, что, как показано в  ~\cite[§7]{FK}, достаточно для доказательства теоремы \ref{uslov}. Теорема доказана.
\section{Функция Гуда и теорема Кусика}\label{good function and kusic th}
Как известно, можно по-другому подойти к определению величины $\delta_{\A}.$ Для этого следует фиксировать алфавит $\A$ вида
\begin{equation}\label{specalpha}
\left\{1,2,\ldots,A-1,A\right\}
\end{equation}
и рассмотреть следующие величины и параметры. Через $V_{\A}(k)$ обозначим множество слов длины $k$
$$
V_{\A}(k)=\left\{(d_1,d_2,\ldots,d_{k})\Bigl| 1\le d_{j}\le A ,\, j=1,\ldots,k\right\},
$$
а через
$V_{\A}=\bigcup_{k\ge1}V_{\A}(k)$
~--множество всех конечных слов. Далее для каждого слова $D=(d_1,d_2,\ldots,d_{k})\in V_{\A}(k)$ берется его континуант ${\langle D\rangle=\langle d_1,d_2,\ldots,d_{k}\rangle}$ (т.е. знаменатель цепной дроби $[d_1,\ldots,d_k]$) и для каждого $s>0$ рассматривается сумма
\begin{equation}\label{Good-zeta-k}
\zeta_k(s,\A)=\sum_{D\in V_{\A}(k)}\langle D\rangle^{-s}.
\end{equation}
Наконец, из этих сумм составляется ряд
\begin{equation}\label{Good-zeta}
\zeta(s,\A)=\sum_{k=1}^{\infty}\zeta_k(s,\A),
\end{equation}
называемый $\zeta$~-- функцией Гуда, сходящийся или расходящийся в зависимости от величины $s.$ Точная нижняя грань тех значений $s,$ при которых ряд \eqref{Good-zeta} сходится, называется абсциссой сходимости этого ряда.\par
Уточняя результаты Гуда~\cite{Good}, Кусик ~(\cite{Cusik1},\,\cite{Cusik2},\,\cite{Cusik3}) доказал теорему о том, что эта абсцисса сходимости равна $2\delta_{\A}$, где, напомним, $\A$ имеет вид \eqref{specalpha}.\par
Хотя для произвольного алфавита $\A$ аналогичное свойство не известно, тем не менее, можно надеяться, что оно когда-нибудь будет доказано. Для нас, однако, абсцисса сходимости ряда \eqref{Good-zeta} важна уже сама по себе, независимо от того, как именно она связана с понятием хаусдорфовой размерности. Дело в том, что следующая ниже теорема \ref{Kan-Hen}, обобщающая результат Хенсли для алфавитов вида \eqref{specalpha}, играет важную роль в построении ансамбля $\Omega_N.$ Отметим, что в ~\cite{BK} аналог теоремы \ref{Kan-Hen} доказан с помощью спектральной теории. Однако наше доказательство теоремы \ref{Kan-Hen} требует определения параметра $\delta_{\A}$  из \eqref{KFcondition} через абсциссу сходимости. Итак, по определению, ряд
\begin{gather}\label{Good-zeta1}
\zeta(s,\A)=\sum_{k=1}^{\infty}\sum_{D\in V_{\A}(k)}\langle D\rangle^{-s}
\end{gather}
\begin{flushleft}
~--при всех\, $s>2\delta_{\A}$~-- сходится,\\
~--при\, $0<s<2\delta_{\A}$~-- расходится,
\end{flushleft}
чем и определяется величина $\delta_{\A}.$  Ситуация, когда $s=2\delta_{\A},$ будет рассмотрена ниже в §\ref{method of hensley} (забегая вперед, скажем, что ряд \eqref{Good-zeta1} расходится и при $s=2\delta_{\A};$ это свойство будет для нас важно).

\section{Метод Хенсли}\label{method of hensley}
Прежде чем выяснять, каково количество континуантов, не превосходящих $N,$ полезно разобраться с вопросом о том, сколь много существует цепных дробей, знаменатель которых не превосходит $N.$ (Несмотря на схожесть вопросов, в ответе на второй из них каждый континуант должен учитываться вместе со своей кратностью.) Для алфавита $\A$ вида \eqref{specalpha} решение этого вопроса содержится в следующей теореме Хенсли ~\cite{Hen2}.
\begin{Th}\label{hensley th}(\textbf{Хенсли ~\cite{Hen2}})
Для действительных $x\ge1,$ для алфавита $\A$ вида \eqref{specalpha} имеет место соотношение
\begin{equation}
\#\left\{D\in V_{\A}\Bigl|\langle D\rangle\le x\right\}\asymp x^{2\delta_{\A}}.
\end{equation}
\end{Th}
Остаток данного параграфа нужен для того, чтобы перейти от алфавита $\A$ вида \eqref{specalpha} к общему случаю. Обобщение этой теоремы Хенсли на случай произвольного алфавита основано на обобщении  его же метода: нужно лишь переписать его доказательство в других обозначениях, иногда упрощая выкладки. К реализации этого плана и приступаем.\par
Рассмотрим произвольный алфавит $\A,$ для которого, как всегда, $|\A|\ge2,$ и вернемся к обозначениям §\ref{good function and kusic th}. Определим для $0\le s\le2, k\in\N$ функцию  $g_{\A}(k,s),$ положив
\begin{equation}\label{def g}
g_{\A}(k,s)=\log\zeta_k(s,\A)-\log2.
\end{equation}
Будем также писать $g(k)$ вместо $g_{\A}(k,s),$ если параметры $s$ и $\A$ ясны из контекста.
\begin{Le}\label{prop g}
Для $s\in[0,2], j,k,r\in\N$ выполнено:
\begin{equation}\label{prop1 g}
|g(j+k)-(g(j)+g(k))|\le\log2,
\end{equation}
\begin{equation}\label{prop2 g}
|g(rj)-rg(j)|\le(r-1)\log2,
\end{equation}
\begin{equation}\label{prop3 g}
|\frac{g(r)}{r}-g(1)|\le\log2.
\end{equation}
\begin{proof}
Ввиду неравенства
\begin{equation}\label{continuant inequality}
\langle D\rangle\langle B\rangle\le
\langle D,B\rangle\le
2\langle D\rangle\langle B\rangle,
\end{equation}
имеем:
\begin{gather}\label{prop1 for1}
2^{-s}\sum_{D\in V_{\A}(j)}\langle D\rangle^{-s}\sum_{B\in V_{\A}(k)}\langle B\rangle^{-s}\le
\sum_{D\in V_{\A}(j)}\sum_{B\in V_{\A}(k)}\langle D,B\rangle^{-s}\le
\sum_{D\in V_{\A}(j)}\langle D\rangle^{-s}\sum_{B\in V_{\A}(k)}\langle B\rangle^{-s},
\end{gather}
или, в терминах функции $\zeta_k(s,\A)$ из \eqref{Good-zeta-k},
\begin{equation}\label{prop1 for2}
2^{-s}\zeta_j(s,\A)\zeta_k(s,\A)\le\zeta_{j+k}(s,\A)\le\zeta_j(s,\A)\zeta_k(s,\A).
\end{equation}
Так как $s\le 2$, то $2^{-s}\ge\frac{1}{4}.$ Поэтому, логарифмируя неравенство \eqref{prop1 for2} и переходя к функции $g(k),$ получаем неравенство \eqref{prop1 g}.
Обобщая его на случай $r\ge2$ слагаемых индукцией по $r,$ получаем:
\begin{equation}\label{prop1 for3}
\left|g\left(\sum_{n=1}^rj_n\right)-\sum_{n=1}^rg(j_n)\right|\le(r-1)\log2
\end{equation}
В частности, при $j_1=\ldots=j_r=j$ получаем неравенство \eqref{prop2 g}. Наконец, подставляя в доказанное неравенство \eqref{prop2 g} значение $j=1$ и производя деление на $r,$ получаем неравенство \eqref{prop3 g}.
\end{proof}
\end{Le}
Определим теперь функции $L(s,\A)$ и $M(s,\A)$ для $s\in[0,2],$ положив:
\begin{equation}\label{def L}
L(s,\A)=\limsup_{r\rightarrow\infty}r^{-1}g(r),
\end{equation}
\begin{equation}\label{def M}
M(s,\A)=\liminf_{r\rightarrow\infty}r^{-1}g(r).
\end{equation}
Корректность определения этих функций уже установлена в неравенстве \eqref{prop3 g}.

\begin{Le}\label{prop L,M}
Функции $L(s,\A)$ и $M(s,\A)$ непрерывны по $s$ на отрезке $[0,2]$ и для всех натуральных $r$ удовлетворяют неравенству
\begin{equation}\label{prop L,g,M}
-\log2+rL(s,\A)\le g(r)\le rM(s,\A)+\log2.
\end{equation}
\begin{proof}
Докажем сначала нижнюю из оценок в \eqref{prop L,g,M}. Для этого предположим противное: пусть для некоторых $j\ge1$и $\epsilon>0$ выполнено неравенство
\begin{equation}\label{prop L,g,M for1}
g(j)\le j(L(s,\A)-\epsilon)-\log2.
\end{equation}
Тогда, ввиду \eqref{prop2 g}, для любого $r\ge1$ получаем:
\begin{equation}\label{prop L,g,M for2}
g(rj)\le rg(j)+(r-1)\log2\le rj(L(s,\A)-\epsilon)-\log2.
\end{equation}
Но любое $m\in\N$ представимо в виде $m=rj+t,\, 0\le t<j.$ Следовательно, ввиду неравенств \eqref{prop1 g} и \eqref{prop L,g,M for2},
\begin{gather}\label{prop L,g,M for3}\notag
g(m)=g(rj+t)\le g(rj)+g(t)+\log2\le
rj(L(s,\A)-\epsilon)+g(t)=\\=
m(L(s,\A)-\epsilon)+g(t)-t(L(s,\A)-\epsilon).
\end{gather}
Так как $t$ принимает лишь конечное число значений, то при $m\rightarrow\infty$ получаем:
\begin{gather}\label{prop L,g,M for3}
g(m)\le m(L(s,\A)-\epsilon'), \,\epsilon'>0,
\end{gather}
что противоречит определению функции $L(s,\A)$ равенством \eqref{def L}, и нижняя оценка в \eqref{prop L,g,M} доказана. Ввиду очевидной симметрии формул, верхняя из оценок в \eqref{prop L,g,M} доказывается полностью аналогично.\par
Используя доказанное неравенство \eqref{prop L,g,M}, получаем, что для любого $r$
$$
M(s,\A)\le L(s,\A)\le M(s,\A)+\frac{2}{r}\log2,
$$
что при $r\rightarrow\infty$ дает:
\begin{equation}\label{L=M}
L(s,\A)=M(s,\A).
\end{equation}
Поэтому достаточно доказать непрерывность только функции $L(s,\A).$ Предположим противное: пусть найдутся $\sigma\in[0,2], \epsilon>0$ и последовательность $s_i\mapsto \sigma$, такая что $0\le s_i\le2$ и
\begin{equation}\label{L-L}
|L(s_i,\A)-L(\sigma,\A)|>\epsilon.
\end{equation}
Ввиду равенства \eqref{L=M}, неравенство \eqref{prop L,g,M} можно записать в форме
\begin{equation}\label{L-g/r}
|r^{-1}g_{\A}(r,s)-L(s,\A)|\le r^{-1}\log2.
\end{equation}
Рассматривая неравенство \eqref{L-g/r} для любого фиксированного $r>\frac{4\log2}{\epsilon}$ при $s=s_i$ для $i=1,2,\ldots$ или $s=\sigma,$ получаем:
\begin{gather}\label{L-g/r 2}
|r^{-1}g_{\A}(r,\sigma)-L(\sigma,\A)|\le \frac{\epsilon}{4},\,
|r^{-1}g_{\A}(r,s_i)-L(s_i,\A)|\le \frac{\epsilon}{4}.
\end{gather}
Теперь из неравенств \eqref{L-L} и \eqref{L-g/r 2} с помощью неравенства треугольника выводится неравенство
\begin{gather}\label{g/r-g/r}
r^{-1}|g_{\A}(r,s_i)-g_{\A}(r,\sigma)|>\frac{\epsilon}{2},
\end{gather}
противоречащее непрерывности $g_{\A}(r,s)$ по $s.$ Этим доказана непрерывность $L(s,\A).$
\end{proof}
\end{Le}

\begin{Le}\label{lambda}
Существует положительная функция $\lambda(s,\A),$ непрерывная по $s$ и строго убывающая при $0\le s\le2,$ такая что:
\begin{equation}\label{lambda1}
\lambda(0,\A)=|\A|,
\end{equation}
\begin{equation}\label{lambda2}
\lambda(2,\A)\le1,
\end{equation}
\begin{equation}\label{lambda3}
\lambda(s,\A)=\lim_{r\rightarrow\infty}(\zeta_r(s,\A))^{\frac{1}{r}},
\end{equation}
\begin{equation}\label{lambda4}
\Bigl|\log\lambda(s,\A)-r^{-1}\log\frac{\zeta_r(s,\A)}{2}\Bigr|\le r^{-1}\log2,
\end{equation}
\begin{equation}\label{lambda5}
\lambda(2\delta_{\A},\A)=1;
\end{equation}
кроме того, выполнено соотношение
\begin{equation}\label{lambda6}
\lim_{s\rightarrow2\delta_{\A}+0}\zeta(s,\A)=+\infty
\end{equation}
\begin{proof}
Положим
\begin{equation}\label{lamb1}
\lambda(s,\A)=\exp(L(s,\A)).
\end{equation}
Тогда положительность $\lambda(s,\A)$ очевидна, а непрерывность следует из непрерывности $L(s,\A),$ доказанной в лемме \ref{prop L,M}. Неравенство \eqref{lambda4} представляет собой записанное в других терминах неравенство \eqref{L-g/r}. Из неравенства \eqref{lambda4} легко получить, что
\begin{equation}\label{lamb2}
\zeta_r(s,\A)\asymp(\lambda(s,\A))^r,
\end{equation}
откуда сразу следует формула \eqref{lambda3}. Из последней тривиально следует \eqref{lambda1}: при $s=0$ все слагаемые в сумме \eqref{Good-zeta-k} равны единице.\par
Далее, из \eqref{lamb2} видно, что ряд $\zeta(s,\A)$ сходится тогда и только тогда, когда ${\lambda(s,\A)<1.}$ Поэтому равенство \eqref{lambda5} выполнено по определению числа $\delta_{\A}$ через абсциссу сходимости. Поскольку, ввиду \eqref{lamb2}, по формуле для суммы геометрической прогрессии
$$\zeta(s,\A)\asymp(1-\lambda(s,\A))^{-1},$$
то неравенство \eqref{lambda6} также доказано.\par
Докажем убывание $\lambda(s,\A)$ по $s.$ Для этого рассмотрим неравенство
\begin{equation}\label{7-29}
\langle D\rangle\ge\left(\frac{\sqrt{5}+1}{2}\right)^{r-1},
\end{equation}
выполненное для всех $D\in V_{\A}(r),$ поскольку
\begin{equation}\label{7-30}
\langle D\rangle\ge\langle \,\underbrace{1,1,\ldots,1}_{r}\,\rangle.
\end{equation}
Поэтому для любого фиксированного $\epsilon>0,$ для всех достаточно больших $r,$ из $D\in V_{\A}(r)$ следует
$$
\langle D\rangle^{-(s+\epsilon)}\le\langle D\rangle^{-s}\left(\frac{\sqrt{5}+1}{2}\right)^{-\epsilon(r-1)}\le\frac{1}{10}
\langle D\rangle^{-s}
$$
Следовательно, ввиду \eqref{Good-zeta-k},
$$
\log\zeta_r(s,\A)-\log\zeta_r(s+\epsilon,\A)\ge\log10.
$$
Ввиду \eqref{lambda4}, применяя неравенство треугольника, получаем
\begin{gather*}
\log\lambda(s,\A)-\log\lambda(s+\epsilon,\A)\ge r^{-1}(\log\zeta_r(s,\A)-\log\zeta_r(s+\epsilon,\A))-\\-
\Bigl|\log\lambda(s,\A)-r^{-1}\log\frac{\zeta_r(s,\A)}{2}\Bigr|-
\Bigl|\log\lambda(s+\epsilon,\A)-r^{-1}\log\frac{\zeta_r(s+\epsilon,\A)}{2}\Bigr|\ge\\\ge
r^{-1}\log10-2r^{-1}\log2>0,
\end{gather*}
что доказывает убывание $\lambda(s,\A)$ по $s.$\par
Осталось доказать неравенство \eqref{lambda2}. Пусть
\begin{gather}\label{maxA}
A=\max\A
\end{gather}
~--наибольший элемент алфавита $\A,$ тогда рассмотрим алфавит $\A'={1,2,\ldots,A}.$ Из сопоставления слагаемых в суммах $\zeta_r(s,\A)$ и $\zeta_r(s,\A')$ очевидно следует неравенство $\delta_{\A}\le\delta_{\A'}.$ Но неравенство $2\delta_{\A'}<2$ было доказано Хенсли \cite[стр.375-377, лемма 1]{Hen2}. Следовательно, $2\delta_{\A}<2,$ и из убывания функции $\lambda(s,\A)$ получаем:
$$
\lambda(2,\A)<\lambda(2\delta_{\A},\A)=1.
$$
\end{proof}
\end{Le}
Договоримся далее рассматривать в алфавите $\A$ слова только четной длины (состоящие из четного числа букв). Можно взглянуть на такие слова чуть иначе, представив себе, что они составлены из букв алфавита $(\A,\A),$ то есть из пар вида $(a,b),$ где $a,b\in\A.$ Обозначим алфавит $(\A,\A)$ через $\A^2,$ а множество слов четной длины в алфавите $\A$~-- через $V_{\A^2}.$ Для алфавита $\A^2$ теперь необходимо ввести ряд объектов, таких как множество конечных цепных дробей $\R_{\A^2}$ (то есть, множество конечных цепных дробей от последовательностей из $V_{\A^2}$), функцию Гуда
\begin{gather}\label{Good-zeta2}
\zeta(s,\A^2)=\sum_{r=1}^{\infty}\sum_{D\in V_{\A}(2r)}\langle D\rangle^{-s}
\end{gather}
и ее абсциссу сходимости $2\delta_{\A^2}.$
\begin{Zam}\label{zam7/1}
Из равенств \eqref{lambda3} и \eqref{lambda5} следует, что $2\delta_{\A^2}=2\delta_{\A},$ так как предел сходящейся последовательности равен пределу по ее подпоследовательности четных индексов. По тем же причинам равенство \eqref{lambda6} остается верным при замене  $\A$ на $\A^2$.
\end{Zam}
Положим также
$$\R_{\A^2}(N)=\left\{D\in V_{\A^2}\Bigl|\langle D\rangle\le N\right\},$$
$$
F_{\A}(x)=\#\left\{D\in V_{\A^2}\Bigl|\langle D\rangle\le x\right\}=\#\R_{\A^2}(x),
$$
и пусть $D_{\A^2}(N)$~--множество знаменателей для дробей из $\R_{\A^2}(N).$
\begin{Th}\label{Kan-Hen}
Пусть $\delta_{\A}>\frac{1}{2},$ тогда для любого $x\ge4A^2$ выполнено:
\begin{gather}\label{F_A less}
\frac{1}{32A^4}x^{2\delta_{\A}}\le F_{\A}(x)-F_{\A}\left(\frac{x}{4A^2}\right)\le F_{\A}(x)\le8 x^{2\delta_{\A}};
\end{gather}
\begin{proof}
Докажем сначала верхнюю оценку в \eqref{F_A less}. Для этого установим оценку
\begin{gather}\label{F_A 1}
F_{\A}(x)-F_{\A}\left(\frac{x}{2}\right)\le4x^{2\delta_{\A}},
\end{gather}
откуда нужное неравенство будет следовать тривиально.\par
Фиксируем действительное $x\ge4A^2$ и рассмотрим множество последовательностей  $D$ из $V_{\A^2},$ таких что $\langle D\rangle>x,$ скажем,
\begin{gather}\label{F_A 2}
D=\{d_1,d_2,\ldots,d_{2r}\}.
\end{gather}
Каждая такая последовательность $D$ может быть единственным образом записана в виде $D=(E,F),$ или
\begin{gather}\label{F_A 3}
\{d_1,d_2,\ldots,d_{2r}\}=\{e_1,e_2,\ldots,e_{2k},f_1,f_2,\ldots,f_{2j}\},
\end{gather}
так чтобы выполнялось неравенство
\begin{gather}\label{F_A 4}
\langle E\rangle\le x<\langle E,f_1,f_2\rangle.
\end{gather}
Оценим континуант $\langle E,f_1,f_2\rangle$ снизу. Ввиду \eqref{continuant inequality}, имеем:
\begin{gather}\label{F_A 7}
\langle E,f_1,f_2\rangle\ge\langle E\rangle\langle f_1,f_2\rangle\ge2\langle E\rangle,
\end{gather}
поскольку
\begin{gather*}
\langle f_1,f_2\rangle=f_1f_2+1\ge2.
\end{gather*}
Поэтому неравенство  \eqref{F_A 4} гарантированно следует из неравенства
$$
\frac{x}{2}<\langle E\rangle\le x.
$$
Следовательно, для любого $s>2\delta_{\A}$ имеем:
\begin{gather}\notag
\zeta(s,\A^2)=\sum_{D\in V_{\A^2}}\langle D\rangle^{-s}\ge\sum_{D\in V_{\A^2}}\langle D\rangle^{-s}\1_{\{\langle D\rangle>x\}}\ge\\\ge
\sum_{F\in V_{\A^2}}\sum_{E\in V_{\A^2}}\langle E,F\rangle^{-s}\1_{
\{\frac{x}{2}<\langle E\rangle\le x\}
},\label{F_A 9}
\end{gather}
где, здесь и далее,
\begin{gather}\label{F_A 10}
\1_{\{S\}}=
\left\{
              \begin{array}{ll}
                1, & \hbox{если условие $S$ выполнено;} \\
                0, & \hbox{иначе.}
              \end{array}
\right.
\end{gather}
Применяя верхнюю из оценок \eqref{continuant inequality} и используя условие $s\le2,$ продолжим цепочку неравенств \eqref{F_A 9}:
\begin{gather*}
\zeta(s,\A^2)\ge
\frac{1}{4}\sum_{F\in V_{\A^2}}\langle F\rangle^{-s}\sum_{E\in V_{\A^2}}\langle E\rangle^{-s}
\1_{
\{\frac{x}{2}<\langle E\rangle\le x\}
}\ge
\frac{1}{4}\zeta(s,\A^2)x^{-s}\left(F_{\A}(x)-F_{\A}\left(\frac{x}{2}\right)\right),
\end{gather*}
откуда, деля обе части на $\zeta(s,\A^2)$ при $s>2\delta_{\A},$ получаем:
$$
F_{\A}(x)-F_{\A}\left(\frac{x}{2}\right)\le4x^{s}
$$
Ввиду произвольности числа $s,$ такого что $s>2\delta_{\A},$ доказано неравенство \eqref{F_A 1}. Применяя неравенство \eqref{F_A 1}, при $\delta_{\A}>\frac{1}{2}$ получаем:
\begin{gather*}
F_{\A}(x)\le\sum_{k\ge0}\left(F_{\A}\left(\frac{x}{2^k}\right)-
F_{\A}\left(\frac{x}{2^{k+1}}\right)\right)\le\\\le
4\sum_{k\ge0}\left(\frac{x}{2^k}\right)^{2\delta_{\A}}\le
4x^{2\delta_{\A}}\frac{1}{1-2^{-1}}\le8x^{2\delta_{\A}},
\end{gather*}
и оценка \eqref{F_A less} доказана.\par
Докажем теперь нижнюю оценку в \eqref{F_A less}.Снова возьмем произвольное $x\ge4A^2$ и рассмотрим те $D\in V_{\A^2},$ для которых $\langle D\rangle>x,$ и пусть $D$ будет, как в \eqref{F_A 2}.
Поскольку, ввиду неравенства \eqref{continuant inequality},
\begin{gather*}
\langle E,f_1,f_2\rangle\le2\langle E\rangle\langle f_1,f_2\rangle<4\langle E\rangle A^2,
\end{gather*}
то неравенство \eqref{F_A 4} можно продолжить влево:
\begin{gather*}
\frac{x}{4A^2}<\langle E\rangle\le x.
\end{gather*}
Учтем также замечание \ref{zam7/1}, согласно которому $\zeta(s,\A^2)\rightarrow+\infty$ при $s\rightarrow2\delta_{\A}+0.$ Следовательно, для такого $s$ имеем:
\begin{gather}\notag
\zeta(s,\A^2)=(1+o(1))\sum_{D\in V_{\A^2}}\langle D\rangle^{-s}\1_{\{\langle D\rangle>x\}}\le\\\le
2\sum_{E\in V_{\A^2}}\langle E\rangle^{-s}\1_{\{\frac{x}{4A^2}<\langle E\rangle\le x\}}
\sum_{F\in V_{\A^2}}\langle F\rangle^{-s},\label{F_A 11}
\end{gather}
где последний множитель~--сумма по $F$~--возник ввиду нижней из оценок \eqref{continuant inequality}, примененной к континуанту $\langle D\rangle:$
$$
\langle D\rangle=\langle E,F\rangle\ge\langle E\rangle\langle F\rangle.
$$
Заменяя $\langle E\rangle^{-1}$ на $\left(\frac{x}{4A^2}\right)^{-1},$ получаем:
\begin{gather*}
1\le
2\left(\frac{x}{4A^2}\right)^{-s}\left(F_{\A}(x)-F_{\A}\left(\frac{x}{4A^2}\right)\right),
\end{gather*}
откуда, ввиду неравенства $s=2\delta_{\A}+0\le2,$
\begin{gather*}
F_{\A}(x)-F_{\A}\left(\frac{x}{4A^2}\right)\ge\frac{1}{2\left(4A^2\right)^2}x^{s}>\frac{1}{32A^4}x^{2\delta_{\A}}.
\end{gather*}
Теорема доказана.
\end{proof}
\end{Th}


\end{document}